\documentclass[preprint,12pt]{elsarticle}

\usepackage{amssymb}
\usepackage{svg}

\usepackage{float}
\usepackage{subfig}
\usepackage{booktabs}
\usepackage{hyperref}
\usepackage{amsmath}
\usepackage{lineno}
\usepackage{geometry}
\hypersetup{
	colorlinks=true,
	linkcolor=cyan,
	filecolor=blue,      
	urlcolor=red,
	citecolor=green,
}
\usepackage[linesnumbered,ruled,vlined]{algorithm2e}
\usepackage{listings}

\lstdefinestyle{Python}{
    language        =   Python, 
    basicstyle      =   \zihao{-5}\ttfamily,
    numberstyle     =   \zihao{-5}\ttfamily,
    keywordstyle    =   \color{blue},
    keywordstyle    =   [2] \color{teal},
    stringstyle     =   \color{magenta},
    commentstyle    =   \color{red}\ttfamily,
    breaklines      =   true,   
    columns         =   fixed,  
    basewidth       =   0.5em,
}

\journal{journal}

\begin{document}

\begin{frontmatter}



\title{A Discrete Neural Operator with Adaptive Sampling for Surrogate Modeling of Parametric Transient Darcy Flows in Porous Media}


\address[1]{Research Center for Mathematics and Interdisciplinary Sciences, Shandong University, Qingdao, Shandong Province, 266237, China}
\address[2]{Frontiers Science Center for Nonlinear Expectations, Minister of Education, Shandong University, Qingdao, Shandong Province, 266237, China}
\address[4]{School of Petroleum Engineering, China University of Petroleum (East China), Qingdao, Shandong Province, 266580, China}
\address[3]{Institute of Mathematical Sciences, ShanghaiTech University, Pudong, Shanghai, 201210, China}
\address[5]{School of Civil Engineering, Qingdao University of Technology, Qingdao, Shandong Province, 266520, China}

\author[1,2]{Zhenglong Chen}
\author[1,2]{Zhao Zhang\corref{cor1}}
\cortext[cor1]{Corresponding author}
\ead{zhaozhang@sdu.edu.cn}
\author[4]{Xia Yan}
\author[3]{Jiayu Zhai}
\author[5]{Piyang Liu}
\author[4,5]{Kai Zhang}


\begin{abstract}
This study proposes a new discrete neural operator for surrogate modeling of transient Darcy flow fields in heterogeneous porous media with random parameters. The new method integrates temporal encoding, operator learning and UNet to approximate the mapping between vector spaces of random parameter and spatiotemporal flow fields. The new discrete neural operator can achieve higher prediction accuracy than the SOTA attention-residual-UNet structure. Derived from the finite volume method, the transmissibility matrices rather than permeability is adopted as the inputs of surrogates to enhance the prediction accuracy further. To increase sampling efficiency, a generative latent space adaptive sampling method is developed employing the Gaussian mixture model for density estimation of generalization error.  Validation is conducted on test cases of 2D/3D single- and two-phase Darcy flow field prediction. Results reveal consistent enhancement in prediction accuracy given limited training set.

\end{abstract}
\begin{keyword}
    Discrete Neural Operator \sep 
    Adaptive Sampling \sep 
    Surrogate modeling \sep 
    Darcy flow \sep
    porous media
\end{keyword}







\end{frontmatter}


\section{Introduction}
\label{intro}
With the advancement of artificial intelligence and deep learning, AI-driven scientific computing has made remarkable progress in solving partial differential equations (PDEs)\cite{RAISSI2019686, article, luDeepXDEDeepLearning2021}. The integration of deep neural networks with PDE solutions has emerged as a transformative paradigm in computational mathematics and numerical analysis.
As fundamental mathematical tools for describing physical systems, PDEs are used to describe processes including fluid dynamics, heat transfer, and electromagnetic field propagation. Classical numerical methods such as finite difference, finite element, and finite volume schemes  have long served as the cornerstone for PDE discretization. Nevertheless, these conventional approaches frequently encounter computational bottlenecks when addressing high-dimensional parameter spaces, irregular domains, and high-performance computing scenarios requiring massive parallelization.

Deep learning establishes a data driven paradigm by leveraging the universal approximation properties for function spaces. This approach demonstrates significant potential in offering computationally efficient solutions to PDE-related problems, particularly those involving nonlinear operators and random parameters.
The deployment of neural architectures including convolutional neural networks (CNNs) for spatial correlations, recurrent architectures (RNNs) for temporal evolution, and attention-based transformers for multi-scale interactions has proven effective in constructing PDE surrogates.
A notable advancement is the development of neural operators\cite{kovachkiNeuralOperatorLearninga, wangLearningSolutionOperator2021}, such as DeepOnet \cite{luLearningNonlinearOperators2021} and Fourier Neural Operator \cite{li2020fourier}. These approaches generalizes traditional numerical solvers by learning mappings between function spaces, enabling efficient resolution of PDEs across varying inputs and parameters in continuous spaces. 
For subsurface flows in porous media, physical parameters are highly heterogeneous. The surrogate modeling of such flows can be naturally regarded as image-to-image learning in the discrete domain. Zhao Zhang\cite{zhangPhysicsinformedDeepConvolutional2022, zhangPhysicsinformedConvolutionalNeural2023} built a discrete physics-informed CNN model for simulating transient two-phase Darcy flows and surrogate modeling\cite{TAKBIRIBORUJENI2020104475, KAZEMI2023105960}.

In this paper, we propose a new discrete operator learning structure with adaptive sampling for the surrogate modeling of subsurface Darcy flows with limited training data. Since subsurface reservoir simulation is typically very expensive, we can only obtain limited number of labeled samples. The main idea of this paper is to enhance prediction accuracy given limited training samples by optimizing the network structure and sampling algorithm. Neural operators have been developed to approximate mappings between continuous function spaces. Based on UNet, a new discrete neural operator is proposed for vector spaces associated with finite volume discretisation to enhance the prediction accuracy. Further, adaptive sampling algorithm is built in the latent space to generate samples of higher quality.

The structure proceeds as follows. Section 2 formulates the PDEs governing multiphase Darcy flows in heterogeneous porous media. The third section presents the embedding algorithm incorporating geological prior knowledge via attentional residual convolutional operator learning. Section 4 introduces adaptive sampling and the AROnet structure for spatial-temporal prediction. Section 5 presents the test cases for validation.

\section{Preliminaries and Problem Settings}
This section establishes the discretization scheme for single-phase Darcy flow and two-phase incompressible flow in porous media with heterogeneous random parameter field as well as the neural operator architecture for building surrogates.

\subsection{Single-Phase Darcy Flow}
The governing equation of single-phase slightly compressible Darcy flow is

\begin{equation}
    \phi c_t\frac{\partial P}{\partial t}=\nabla\cdot(\frac{K}{\mu}\nabla P)+f,
    \label{eq:single-phase}
\end{equation}
with initial condition $P(X, t=0) = P_{init}$ for all grid cells, and no-flow condition at all boundaries. \autoref{tab:params} shows the physical meanings of other parameters.

\begin{table}[H]
    \caption{physical parameters}
    \centering
    \begin{tabular}{c c}  
    \toprule
        parameters & physical meaning \\
    \midrule
        $X$ & spatial coordinate \\
        $P$ & pressure\\ 
        $K$ & permeability\\
        $c_t$ & total compressibility \\
        $\mu$ & viscosity \\
        $\phi$ & porosity \\
        $f$ &  source or sink term \\
        $P_{wf}$   &  bottom-hole pressure \\
        $q$     &   volumetric flow rate \\
        
    \bottomrule
    \end{tabular}
    \label{tab:params}
\end{table}

Using the finite volume method(FVM\cite{chen2006computational}), \autoref{eq:single-phase} can be rewritten as
\begin{equation}
    \label{eq:fvm}
    V_{i}\phi C_{i}\frac{P_{i}^{n+1}-P_{i}^{n}}{\Delta t}=\sum_{j}T_{ij}(P_{j}^{n+1}-P_{i}^{n+1})+V_{i}f_{i}^{n+1}\:,
\end{equation}
where, $i, j$ are cell indeces, and $n$ is time step. $T_{ij}$ is the transmissibility between  cell $i$ and $j$ approximated by two-point flux approximation (TPFA) as
\begin{equation}
    \label{eq:Tij}
    T_{ij}=(T_i^{-1}+T_j^{-1})^{-1}\:,
\end{equation}
and $T_i$ is the transmissibility inside cell $i$ towards cell $j$ calculated as 
\begin{equation}
    T_i=\frac{K_iA_{ij}}{\mu d}\:,
\end{equation}
where $A_{ij}$ is the area of boundary face between cell $i$ and cell $j$. The volumetric flow rate $q_{i}^{n+1}=V_{i}f_{i}^{n+1}$ as the source term for cell $i$ containing producing wells\cite{peaceman1978interpretation, brooks1965hydraulic} is modeled as
\begin{equation}
    q_{i}^{n+1}=\mathrm{PI}*(P_{i}^{n+1}-P_{wf})
\end{equation}
where PI is the production index.

\subsection{Two-Phase Darcy Flow}
The governing equations for two-phase slightly compressible Darcy flow, neglecting gravity and capillary pressure \cite{chen2006computational}, are formulated as
\begin{equation}
    \label{eq:2phase}
    \phi\left[S_{\alpha}(c_{r}+c_{\alpha})\frac{\partial P}{\partial t}+\frac{\partial S_{\alpha}}{\partial t}\right]=\nabla\left(\frac{k_{r\alpha}K}{\mu_{\alpha}}\nabla P_{\alpha}\right) + q_{\alpha}\:,
\end{equation}
in which $\alpha=o, w$ represents the non-wetting and wetting phases, respectively. $S_\alpha$ is the phase saturation, $k_{r\alpha}$ is phase relative permeability and $c_\alpha, c_r$ are phase compressibility and rock compressibility. Additionally, $S_w + S_o =1$. The governing equations of two phases are added to obtain 

\begin{equation}
\label{eq:gov}
    \phi(c_r+S_oc_o+S_wc_w)\frac{\partial P}{\partial t}=\nabla\left(\left(\frac{k_{ro}}{\mu_o}+\frac{k_{rw}}{\mu_w}\right)K\nabla P\right) + q_o+q_w\:,
\end{equation}
Using FVM and implicit-pressure explicit-saturation time integration scheme, we have the implicit equation for pressure
\begin{equation}
    \label{eq:fvm-press}
    \frac{V\phi c_t}{\Delta t}(P_i^{n+1}-P_i^n)=\sum_j\lambda_{ij}T_{ij}(P_j^{n+1}-P_i^{n+1})+qV\:,
\end{equation}
where $q=q_{o}+q_{w}$, $\lambda_{ij} = (\lambda_{o}+\lambda_{w})_{ij}=\left(\frac{k_{ro}}{\mu_{o}}+\frac{k_{rw}}{\mu_{w}}\right)_{ij}$ represent total flow rate and total mobility, respectively. Then saturation is updated explicitly as
\begin{equation}
    \label{eq:saturation}
    \frac{V\phi S_{w,i}^{n}C_{w}}{\Delta t}(P_{i}^{n+1}-P_{i}^{n})+V\phi\frac{S_{w,i}^{n+1}-S_{w,i}^{n}}{\Delta t}=\sum_{j}\lambda_{w,ij}T_{ij}(P_{j}^{n+1}-P_{i}^{n+1})+q_{w}\:.
\end{equation}
The relative permeability functions are specified as follows under zero residual oil saturation conditions
\begin{equation}
\begin{aligned}
&k_{rw}=\left(\frac{S_{w}-S_{iw}}{1-S_{iw}}\right)^{4}\:,\\
&k_{ro}=\left(\frac{1-S_{w}}{1-S_{iw}}\right)^{2}\left(1-\left(\frac{S_{w}-S_{iw}}{1-S_{iw}}\right)^{2}\right)~.
\end{aligned}
\end{equation}

\section{Learning PDE Solutions by Neural Network}
Conventional numerical methods for solving PDEs associated with subsurface flows are well established. For parameterized PDEs, numerical simulation is needed for each realization of random parameter fields. The objective of this work is to construct an accurate spatial-temporal surrogate for real-time responses given realizations of heterogeneous parameter fields. The model accepts realizations of random parameters and returns the corresponding PDE solution.

\subsection{Attentional Residual U-net}
The Attentional U-net\cite{oktay2018attention} is a convolutional neural network (CNN) architecture that has been widely used in image segmentation tasks. It consists of an encoder-decoder structure\cite{cho2014learning} with skip connections, allowing for the preservation of spatial information while capturing multi-scale features. The attention mechanism enhances the model's ability to focus on relevant regions in the input data, improving performance in tasks such as semantic segmentation and image generation\cite{ho2020denoising, song2020score}.

Based on U-net, Resnet \cite{he2016deep} and attention mechanism, the Attentional Residual U-net (ARUnet) has been the SOTA for image-to-image learning tasks \cite{huang2021gcaunet,ding2024novel}, which is essentially consistent with approximating the mapping between discretised parameter and solution fields for PDEs.
The residual connect can avoid network degradation which makes the deep network perform worse than the shallow network.

\begin{figure}[htbp]
    \centering
    \subfloat[]{\includegraphics[width=0.99\linewidth]{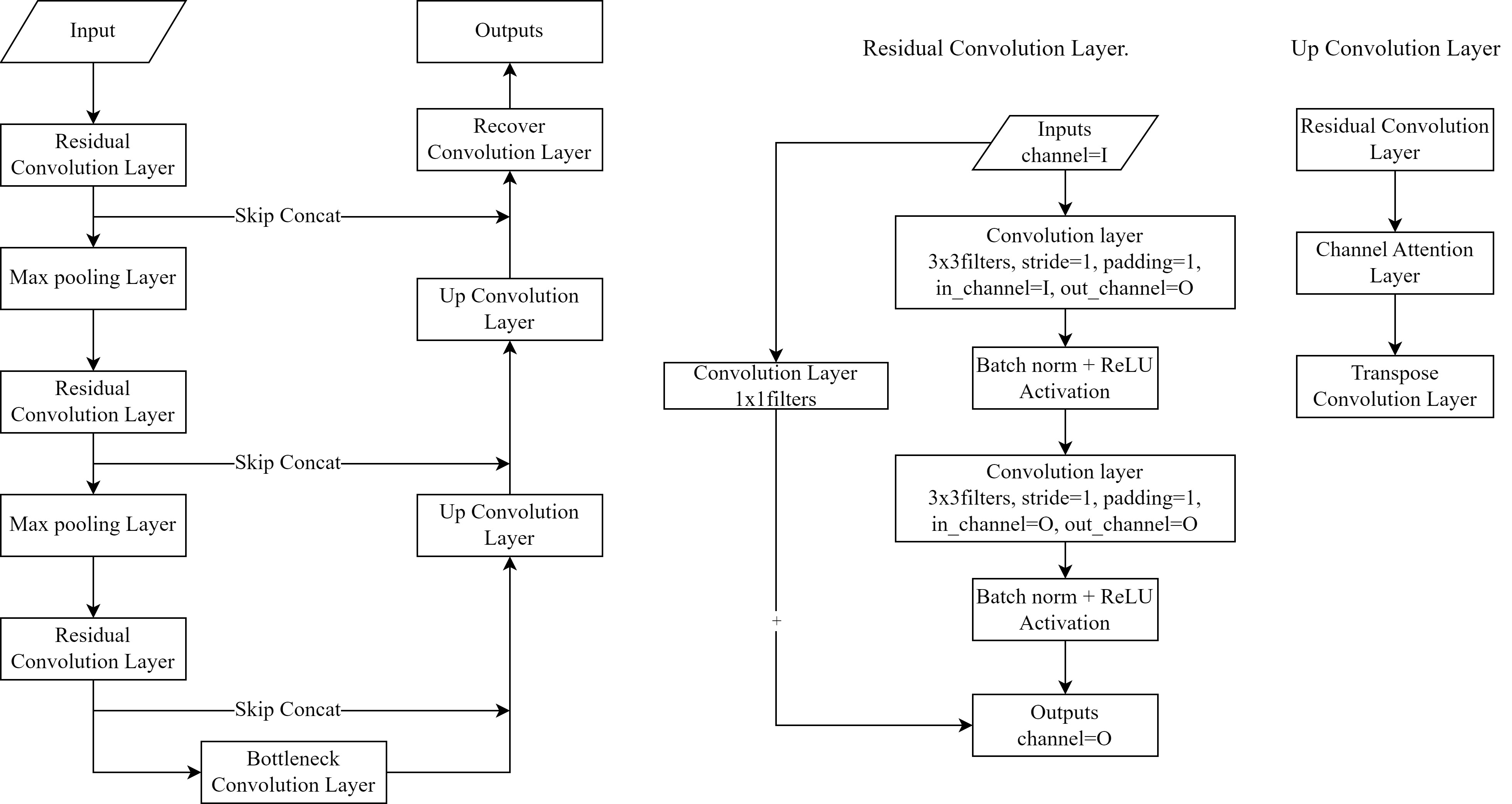}} \\
    
    \subfloat[]{\includegraphics[width=0.99\linewidth]{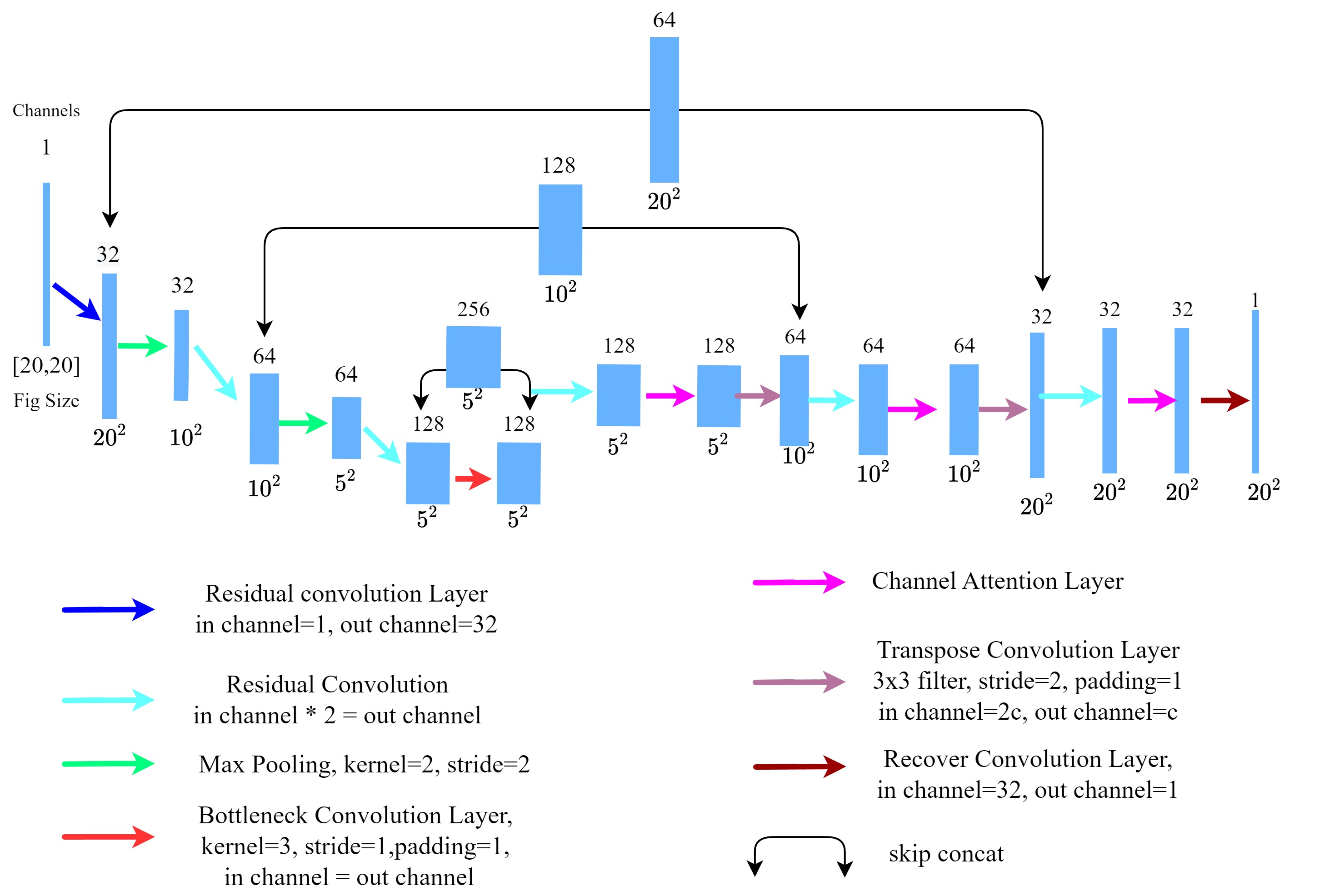}}
    \caption{(a): Attentional Residual U-net structure. (b): An example computational process of ARUnet, using an image data with shape [batch size, channels, hight, weight] = [1, 1, 20, 20] as input of NN.}
    \label{fig:Unet}
\end{figure}

\subsection{Attention Residual Operator}\label{sec:time-embedding}
Based on ARUnet and inspired by deepOnet \cite{luLearningNonlinearOperators2021}, a discrete neural operator, Attention Residual Operator net (AROnet), is developed for operator learning between discrete vector spaces.
AROnet reformulates flow field prediction as an operator learning problem, where the objective is to train a neural operator $G_{\theta}$ that approximates the underlying PDE operator $G$ as follows.
\begin{equation}
    G(u(x))(t) = p(t) , \quad x \in X, t \in [0, T] \:,
\end{equation}
where $u(x)$ represents PDE parameters (e.g., permeability/transmissibility fields) over spatial coordinates $x$, and $t$ denotes time within the solution domain. The objective is to develop a surrogate model $G_{\theta}$ that effectively replaces the operator $G$, allowing for efficient predictions at any given $t$.

For predicting transient flow fields, CNN-LSTM architectures sequentially predict predefined chronological steps. such approaches inherently restrict temporal flexibility and faces difficulties as error accumulates with increasing time steps. 
AROnet addresses these limitations through an operator learning framework using time embedding.  By implementing a sine-cosine scheme, the model disentangles time representation from rigid chronological sequences, enabling simultaneous predictions at any time instance within the solution domain.
Temporal information $t$ is embedded into neural networks by sine-cosine embedding, 
\begin{equation}
\begin{aligned}
    TE_{(t,2i)} &=\sin{(t/10000^{2i/d})}\\
    TE_{(t,2i+1)} &=\cos{(t/10000^{2i/d})}
    \label{eq:TE}
\end{aligned}
\end{equation}
in which $d$ represents the embedding dimension. Here we set $d$ equal to the size of input $u$ to be consistent with the convolutional layers in the NN structure.
The AROnet approximates $G$ as follows.
\begin{equation}
    G_{\theta}(u)\left(t\right)=f(\sum_{k=1}^{q}\underbrace{\mathbf{b}_{k}(u(x_{1}),u(x_{2}),\ldots,u(x_{m}))}_{\mathrm{branch}} * \underbrace{\mathbf{t}_{k}(t)}_{\mathrm{trunk}}) \:,
    \label{eq:Onet}
\end{equation}
where $f$ is a CNN with sigmoid activation layer, while \textbf{branch} and \textbf{trunk} are two different sub-networks. To let $G_{\theta}(u) \approx G(u)$, the loss function is written as 
\begin{equation}
    \mathcal{L}_{Operator}(\theta) = \frac{1}{NM} \sum_{i=1}^N \sum_{j=1}^M \left| G_{\theta}(u^{(i)})(y_j^{(i)}) - G(u^{(i)})(y_j^{(i)}) \right|^2 \:,
\end{equation}
where $N$ represents the number of samples and $M$ for steps of time.

\begin{figure}[htb]
    \centering
    \subfloat[]{\includegraphics[width=0.75\linewidth]{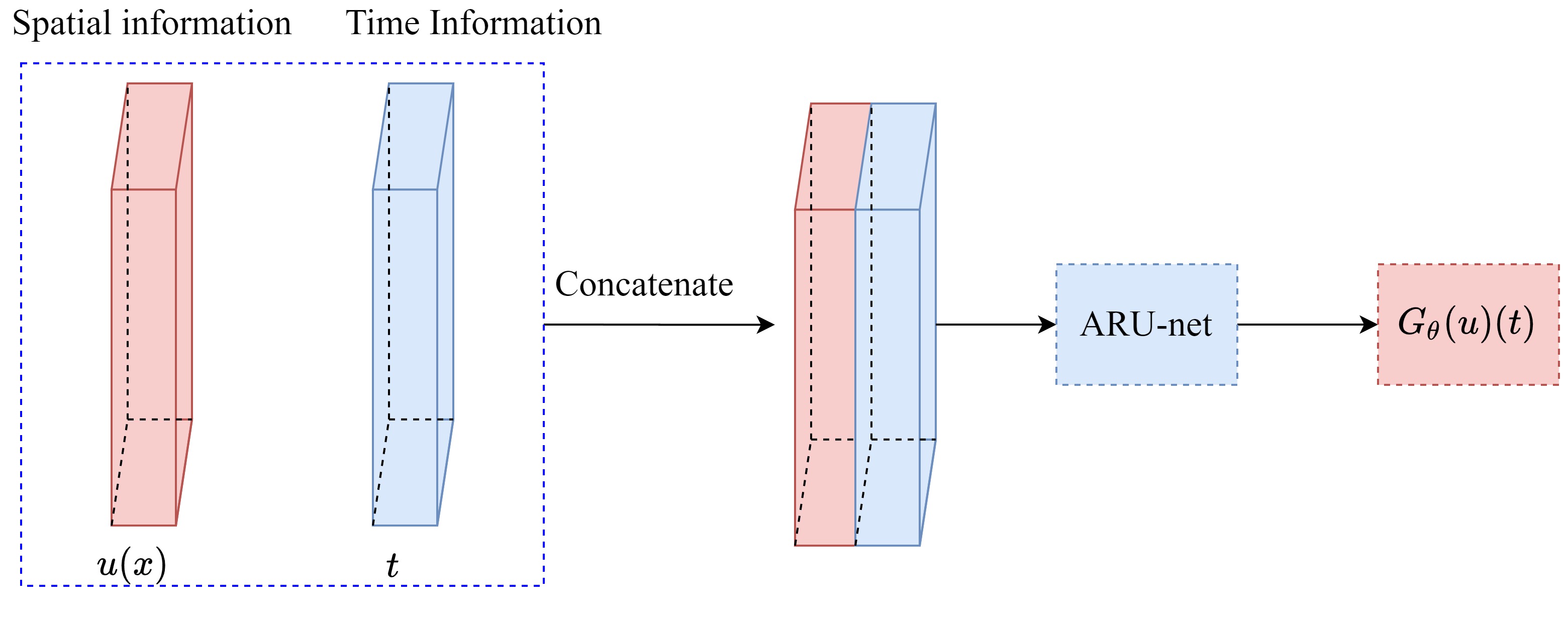}}\\
    \subfloat[]{\includegraphics[width=0.75\linewidth]{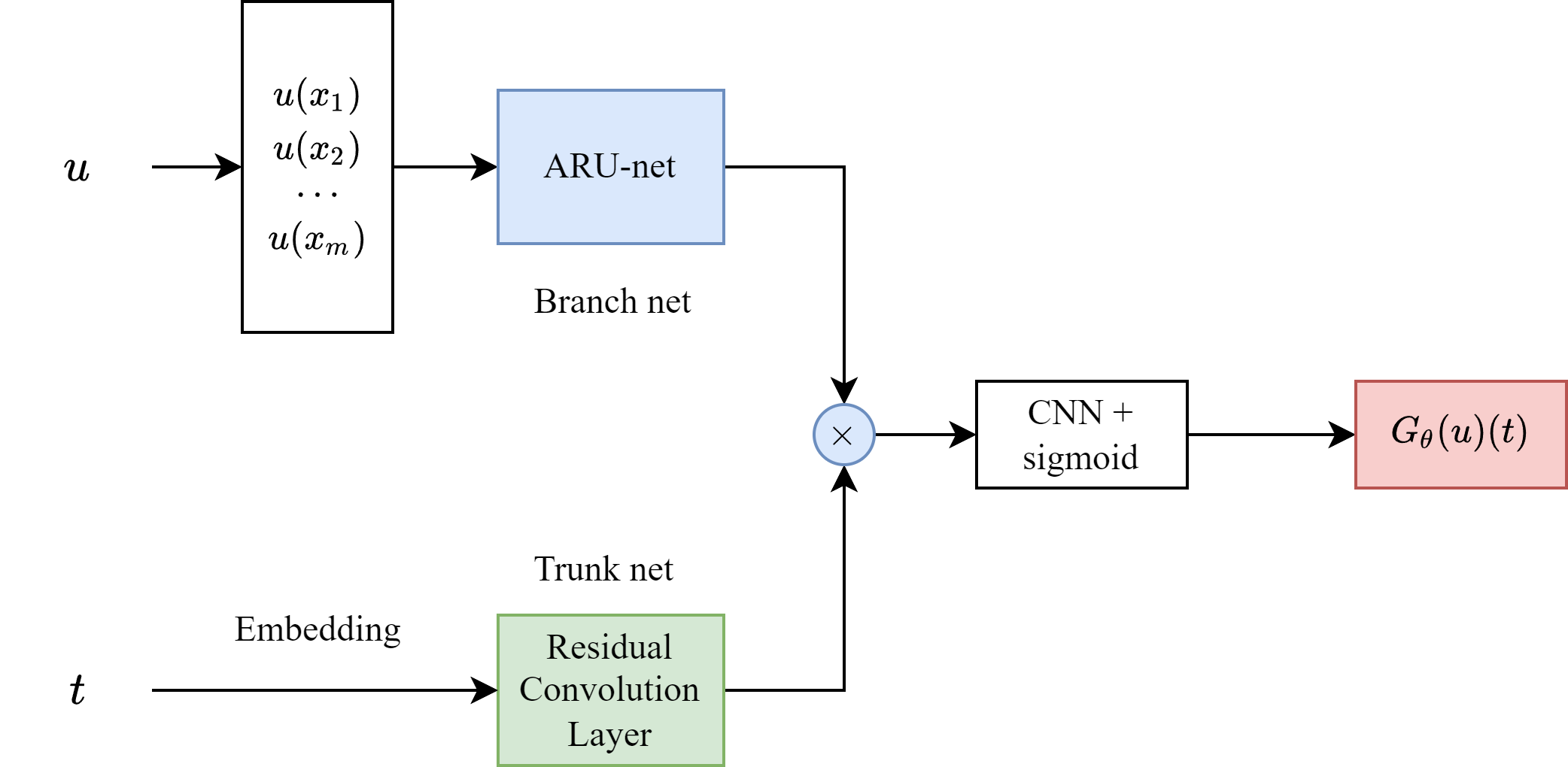}}
    
    \caption{(a): ARUnet Temporal Prediction Framework. (b): AROnet Operator Learning Framework.}
    \label{fig:AROnet}
\end{figure}

As shown in \autoref{fig:AROnet}, both ARUnet and AROnet take the same input and output. The difference is that ARUnet uses a CNN structure to predict the pressure field at each time step, while AROnet learns the operator map $G(u, t) \rightarrow p(t)$ directly where $p(t)$ denotes the labels for outputs.
In ARUnet, parameter $u$ is encoded via convolutional blocks to capture spatial patterns, while discrete time index $t$ is embedded as a positional vector through sine-cosine embedding. Spatial and temporal embeddings are concatenated as different channels. 
AROnet establishes spatiotemporal mapping based on the neural operator architecture. The branch net generates parameter-dependent feature maps from $u$, and the trunk net transforms time-encoded tensor into spatial modulation weights. Branch features and Trunk weights are combined via channel-wise multiplication, followed by depth-separable convolution and ReLU activation to output predictions. The inputs and outputs are both normalized for enhanced prediction accuracy.

The performance of AROnet and ARUnet is compared in \autoref{sec:2D-val}. The experiments results show that the mean relative error and well-block relative error of AROnet is lower than ARUnet.

\subsection{Inputs and Outputs of Neural Network}
The inputs of NN are a set of realizations for the random parameter field of PDEs. In prior studies \cite{zhangPhysicsinformedConvolutionalNeural2023}, the inputs of surrogate models are directly the random parameter field using permeability matrices as an example, and the outputs are corresponding PDE solutions using pressure as an example. While permeability matrices provide straightforward input representations, the influence of transmissibility tensor on solution is more direct according to discretization schemes \autoref{eq:fvm} and \autoref{eq:fvm-press}.

In \autoref{sec:2D-val}, we compares the relative error distributions of predicted pressure field using the same neural operator architecture but with different inputs: one taking permeability tensors $K_{ij}$ as inputs, and the other utilizing transmissibility matrices $T_{ij}$ derived from finite volume discretization. The experiments results shows that, the latter converges faster during training process and achieves higher prediction accuracy and the detail results are shown in \autoref{sec:2D-val}.




\section{Adaptive Sampling in Latent Space}
According to the theoretical framework established in \cite{scott2015multivariate, tangDASPINNsDeepAdaptive2023}, the total error of NN-based surrogates comprises of two parts, i.e., the approximation error arising from model capacity limitations, and statistical error originating from limited training data.
Consider an abstract problem formulation,
\begin{equation}
\mathcal{L}u(x) = s(x), \quad \forall x \in \Omega,
\end{equation}
where $\mathcal{L}$ denotes the governing differential operator containing boundary conditions, $\Omega$ represents the complete sample space, $u(x)$denotes the PDE solution field, and $s(x)$ is the source term. The primary objective of surrogate modeling is to construct a neural network approximation $\mathcal{L}_{\Theta} u(x) = u(x, \Theta)$ parameterized by $\Theta$ to minimize the loss functional
\begin{equation}
    \label{eq:approximate}
    {J}(u(x;\Theta)) = |u(x,\Theta) - s(x)|_{\Omega}^2~.
\end{equation}

Setting $u(\cdot,\Theta_N)$ as the optimal approximation in the training sample space $\Omega_N$, $u(\cdot,\Theta)$ as the optimal approximation in the full sample space $\Omega$, 

\begin{align}
    u(\cdot;\Theta ) &=\arg\min_{x \in \Omega, \Theta}J(u(x;\Theta)),\\
    u(\cdot;\Theta_{N} ) &=\arg\min_{x \in \Omega_N, \Theta}J_{N}(u(x;\Theta)).
\end{align}

Applying expectation operator and Minkowski inequality to the loss functional yields the decomposition as
\begin{equation}
\label{eq:error decomposition}
\mathbb{E}[J(u(x; \Theta))] \leq \underbrace{\mathbb{E}\left(\left \|u(\cdot;\Theta_N) - u(\cdot;\Theta) \right \|_{\Omega}\right)}_{\text{statistic error}} + \underbrace{\mathbb{E}\left(\left \|u(\cdot;\Theta) - s(\cdot)\right \|_{\Omega}\right)}_{\text{approximation error}}~,
\end{equation}
where the approximation error is due to neural network capability, and the
statistic error results from finite training samples. This decomposition reveals two fundamental improvement pathways, i.e.
\begin{itemize}
    \item Enhancing approximation capability by modifying neural network structures.
    \item Optimizing sample distributions via adaptive sampling strategies.
\end{itemize}
The first approach has been extensively studied in literature \cite{luLearningNonlinearOperators2021,zhangPhysicsinformedConvolutionalNeural2023}, where various NN structures have been proposed to enhance the approximation capabilities. The neural operator structure in the current work is developed in \autoref{sec:time-embedding}. In \autoref{sec:opt-sampling}, we focus on the second aspect to develop adaptive sampling strategies that leverage the predicted residual distribution to improve sample distributions.

\subsection{The Objective of Adaptive Sampling} \label{sec:opt-sampling}

The adaptive sampling methodology is formulated for PINNs in \cite{luDeepXDEDeepLearning2021} through their residual-based adaptive refinement (RAR) algorithm, which allocates collocation points in regions of elevated PDE residual. Developments by \cite{tangDeepDensityEstimation2020} established a generative KRnet architecture leveraging normalizing flow models \cite{dinh2016density} for adaptive sampling \cite{tangDASPINNsDeepAdaptive2023}.

The RAR algorithm implements a hard thresholds measure characterized by the discontinuous density, and the $K$ points with the largest residuals will be re-added to the training set, which means that its sampling density $q^*(x)$ is defined as:
\begin{equation}
    \label{eq:rag}
    q^*(x) = \begin{cases}
            \frac{1}{K}, \quad r(x) \geq R ;\\
            0, \quad r(x) < R .
            \end{cases}
\end{equation}
where $R$ is the residual thresholds satisfying $\mu(\{x\in\Omega:r(x)\geq R\})=\frac{K}{N}$. This implies that the sampling density is not continuous, and samples with different residuals are allocated the same probability density, which leads to inaccurate approximation of the true density.
The samples are generally associated with different residual magnitudes. The basic idea of adaptive sampling is to minimize $\mathbb{E}\left( J\left(u(\cdot;\Theta)\right) \right)$ by assigning appropriate sampling density in regions with different prediction errors.

The other track is to use generative model for sample generation. Under the theoretical premise that the residual $r^2(x, \Theta)$ in \autoref{eq:approximate} admits a probability density $p_{\alpha}(x)$ over $\Omega$ that can be chosen. So the key for adaptive sampling is the modeling of $p_{\alpha}(x)$. The optimal design is to have a similar profile with the residual distribution, that is, 
adaptive methods resolve the measure duality between the intrinsic residual distribution $p_{\alpha}(x) \propto r^{2}(x;\Theta)$ and the sampling measure $q(x)$ through parametric density estimation:
\begin{equation}
\min_{q \in \mathcal{P}(\Omega)} D_{\text{KL}}\left( q(x) || \frac{r^2(x;\Theta)p_{\alpha}(x)}{\int_\Omega r^2(y;\Theta)p_{\alpha}(y)dy} \right),
\end{equation}
where $\mathcal{P}(\Omega)$ is the space of probability measures on $\Omega$. The solution $q^*(x) = \arg \min D_{KL}(\cdot)$ recovers the theoretically optimal sampling density. 

To achieve this goal in an automatic process, we use the following minimax scheme \cite{jiao2024gas, tang2023adversarial} to obtain it dynamically with the training process. So loss has two parts for two models. One is to minimize the residuals of the surrogate model on the training samples, which is the optimization parameter $\Theta$. The other is to maximize the residuals for the adaptive sample generater, which is the optimization parameter $\alpha$.

\begin{equation}
\label{eq:EJ}
\begin{aligned}
    \min_{\Theta} \max_{\alpha} \mathbb{E}\left( J\left(u(\cdot;\Theta)\right) \right) = &\int_{\Omega} r^2(x; \Theta) p_{\alpha}(x) \mathrm{d}x \\
    s.t. &\int_{\Omega} p_{\alpha}(x) \mathrm{d}x = 1.
\end{aligned}
\end{equation}
The probability distribution $p_\alpha$ can be parametrized in a standard way \cite{jiao2024gas} using a Gaussian mixtrue model, or using a neural network based model, e.g., normalizing flow \cite{tang2023adversarial}. We use a Gaussian mixture model in this work. The convergence of this minimax framework and its adaptivity property is proved in \cite{tang2023adversarial} using optimal transport theory.


\subsection{Adaptive Sampling in Latent Space via Gaussian Mixture Model}
The Gaussian Mixture Model (GMM) is adopted for density estimation and generative adaptive sampling in the latent space. Two key points for surrogate modeling are addressed as follows.
\begin{enumerate}[(1)]
    \item Adaptive sample generation in regions of higher prediction uncertainty, achieved by iteratively updating the GMM to prioritize underperformed regions of the current surrogate model which is optimized iteratively along with the new samples.
    \item Curse of dimensionality mitigation via latent space sampling. The high-dimensional physical fields are projected into a compressed latent representation using principal component analysis (PCA) as an encoder.
\end{enumerate}
 
\begin{figure}
    \centering
    \includegraphics[width=0.8\linewidth]{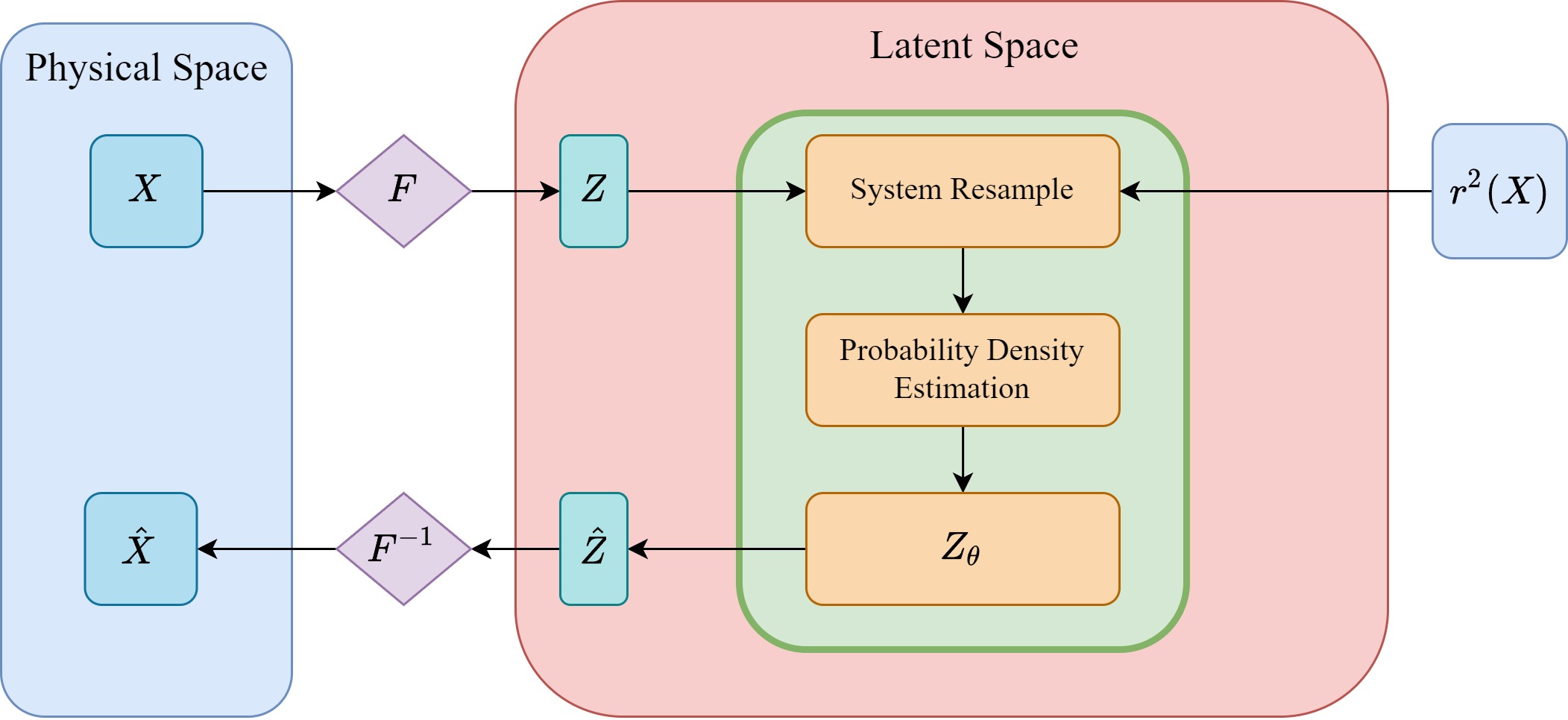}
    \caption{Adaptive Sampling in Latent Space.}
    \label{fig:Latent-AS}
\end{figure}

 Subsequent GMM density estimation and adaptive sampling are executed in the low-dimensional latent space. The explicit density modeling in latent space coupled with decoder-based sample reconstruction is for the balance between computational efficiency and sampling quality in the iterative surrogate optimization process. Adaptive sampling process is shown in \autoref{fig:Latent-AS}. For the test case, the original 1600-dimensional training data is projected into a 70-dimensional latent space ($X \stackrel{F}{\rightarrow} Z, X \in \mathbb{R}^{1600}, Z \in \mathbb{R}^{70}$) using PCA, preserving 95\% cumulative variance:
$$\frac{\sum_{i=1} ^ {70} \lambda _ { i } } { \sum _ { j = 1 } ^ { 100} \lambda _ {j} } = 0 .95,$$
where $\lambda_i$ is the sorted $i$-th eigenvalue of $X^TX$  in descending order.

\begin{enumerate}
    \item The dataset is partitioned into training data  $X_{t}, y_t$ and testing data $X_v, y_v$. The procedure is initiated by training a neural network using an initial training subset $X_{t}^{(0)}, X_{t}^{(0)}\subset X_{t}$, Subsequently, the network generates predicted outputs $\hat{y}_t^{(0)}$ for $X_{t}^{(0)}$, from which residual vectors are computed as
    \begin{equation}
        \label{eq:residual1}
       r^2(X_t) \Leftrightarrow   R_t^{(0)}=\left\|y_t^{(0)}-\hat{y}_t^{(0)}\right\|_2^2,
    \end{equation}
    Each residual component $R_t^{(0),i}$ quantifies the mean squared error (MSE) across all spatial discretization points (mesh grids) associated with sample $i$.
    \begin{equation}
        \label{eq:residual2}
        R_t^{(0),i}=\frac{1}{N_g}\sum_{j=1}^{N_g}\left(\hat{y}_{t,j}^{(0),i}-y_{t,j}^{(0),i}\right)^2
    \end{equation}
    where $N_g$ represents the number of mesh nodes per sample.
    \item The Latent feature set $Z_t^{(0)}$ is obtained by applying PCA to $X_t^{(0)}$, i.e. 
    $$Z_t^{(0)}=\mathbf{W}_{\mathrm{PCA}}^TX_t^{(0)},$$ 
    where $\mathbf{W}_{\mathrm{PCA}}$ denotes the PCA projection matrix. $R_t^{(0)}$ serves as an importance sampling distribution over $Z_t^{(0)}$ and using Systematic resampling to sample $Z_t^{(0)}$ to obtain new samples $\hat{Z}$. Systematic resampling is then performed as \autoref{alg:sys-res}.
    
    \item The resampled data $\hat{Z}$ is used to train a Gaussian Mixture Model (GMM). The GMM is parameterized by $m$ components, 
    \begin{equation}
        \label{eq:gmm}
        p(x)=\sum_{i=1}^{m}w_i\mathcal{N}(x|\mu_i,\Sigma_i),
    \end{equation}
    where $ m \in \mathcal{H}$ is the number of Gaussian distributions in the mixture which is a hyperparameter set. The GMM is trained using the Expectation-Maximization (EM) algorithm, which iteratively estimates the parameters of the GMM until convergence. The optimal component count $m^*$ is determined by minimizing the Bayesian Information Criterion(BIC) over candidate GMMs ${\mathcal{M}_m | m \in \mathcal{H}}$
\begin{equation}
    \label{eq:m}
    m^*=\arg\min_{m\in\mathcal{H}}\mathrm{BIC}(\mathcal{M}_m)
\end{equation}
where the BIC is computed as:
\begin{equation}
    \label{eq:bic}
    \mathrm{BIC}(\mathcal{M}_m)=-2\mathcal{L}_m+d_m\log N
\end{equation}
where $\mathcal{L}_m=\sum_{i=1}^N\log p(\hat{\mathbf{Z}}_i;\Theta_m)$ is the log-likelihood of resampled data $\hat{Z}$, $d_m=m(1+\frac{n(n+1)}2)-1$ is degrees of freedom for n-dimensional GMM with $m$ components and $N$ is the number of resampled points in $\hat{Z}$.

    \item Generate $N$ new samples from BIC-optimized GMM $\mathcal{M}_{m^*}$, denoted as $\hat{Z}$. The new samples are then projected back to the original feature space using the inverse PCA transformation:
    $$\hat{X}=\mathbf{W}_{\mathrm{PCA}} \hat{Z} + \mu_{X}$$
    where $\mathbf{W}_{\mathrm{PCA}}$ is the PCA projection matrix. Labeled data $\hat{y}$ for new samples is computed by solving the PDEs using FVM. The detailed process is shown in \autoref{alg:prob-est}.
\end{enumerate}

\section{Numerical Validation}\label{sec:val}
This Section mainly displays the experimental results after the relevant experiments are carried out using the methods introduced in the previous sections, including the error distribution charts of different model methods, the tables drawn by the relative errors of different models, and the prediction results of test examples.

\subsection{2D Single-Phase Darcy Flow}\label{sec:2D-val}

\begin{figure}[H]
    \centering
    \subfloat[]{\includegraphics[width=0.5\textwidth]{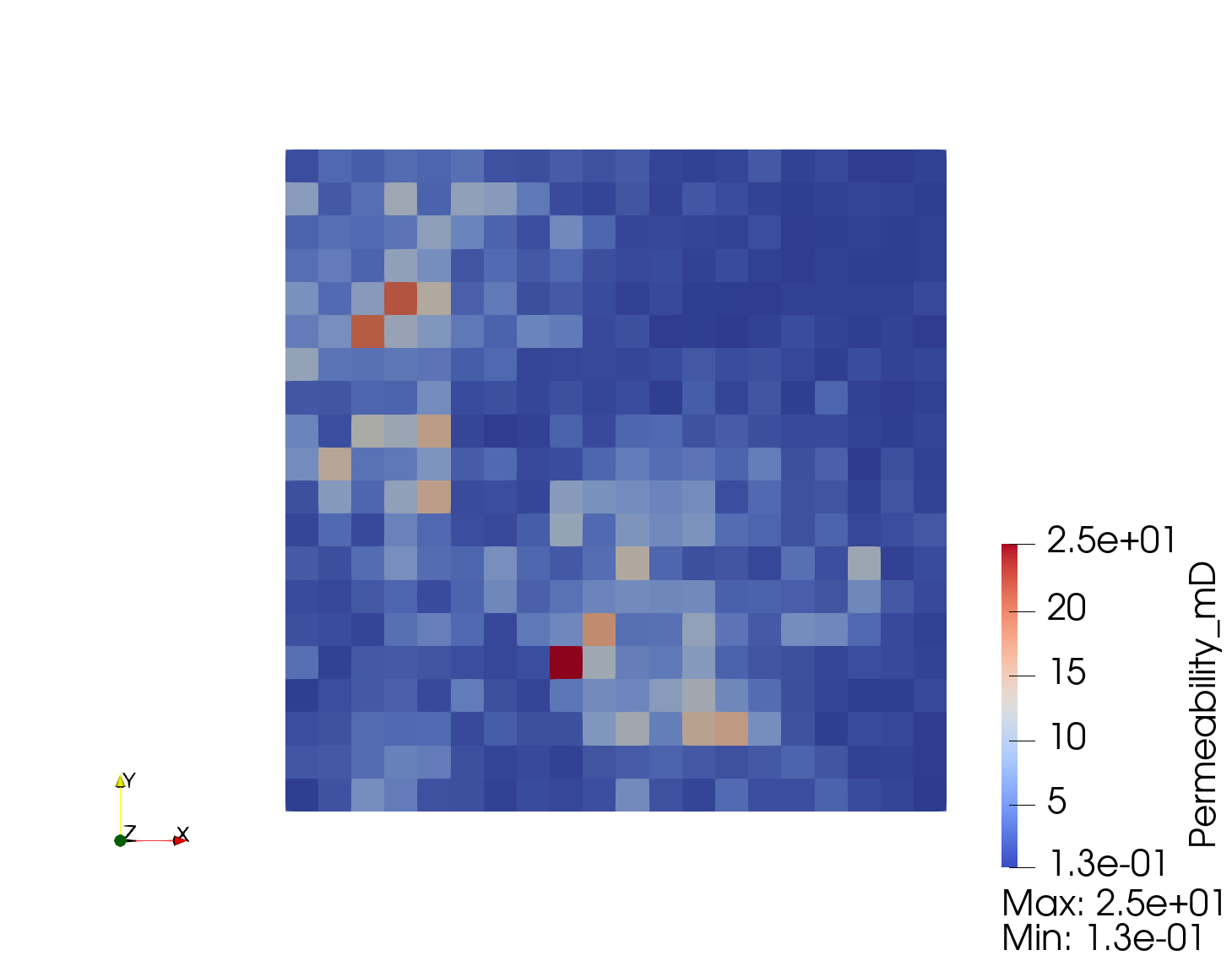}}
    \subfloat[]{\includegraphics[width=0.5\textwidth]{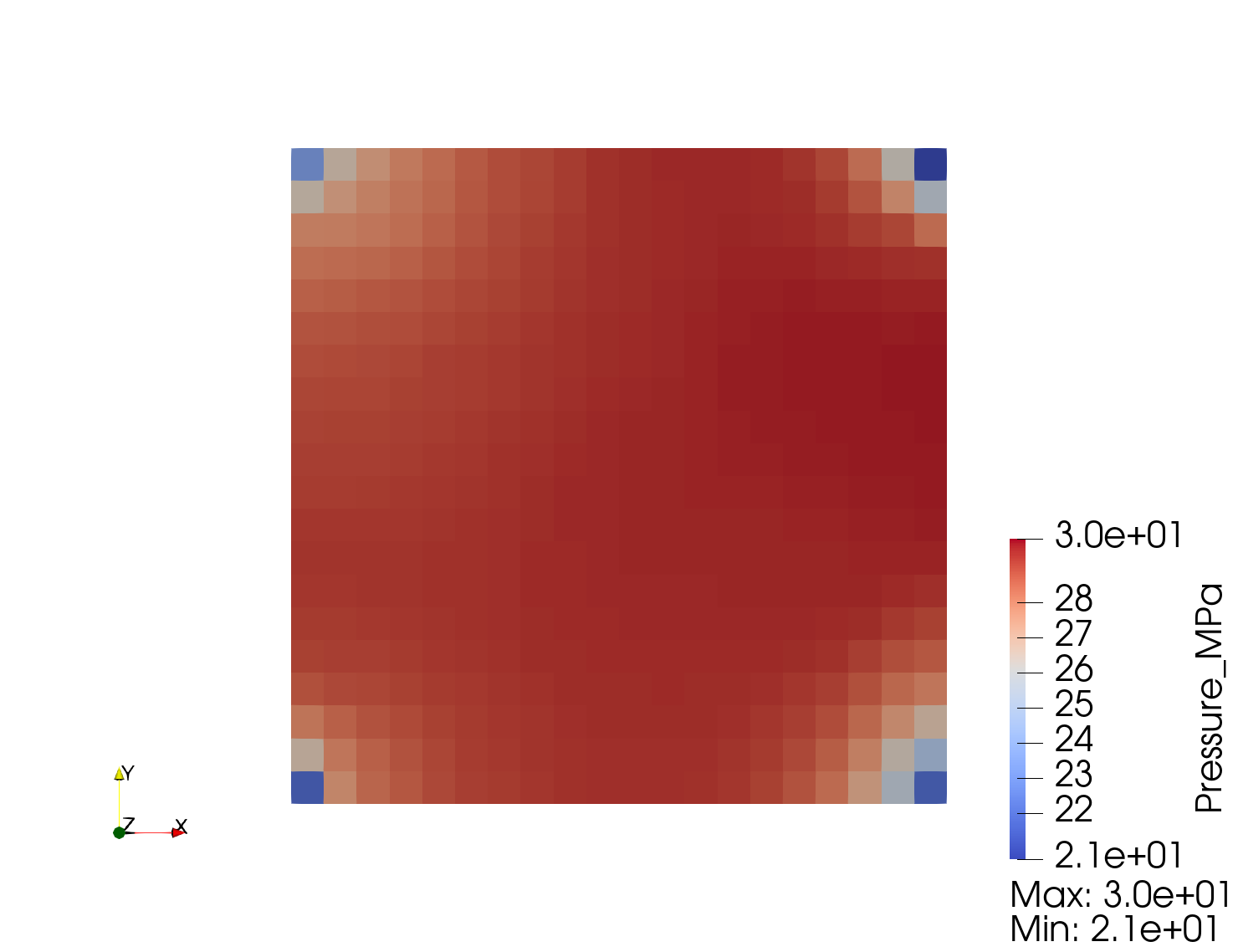}}
    \caption{(a): Heterogeneous permeability field generated by sequential Gaussian simulation. (b): 
    Corresponding pressure computed by FVM at last time step.}
    \label{fig:perm-press}
\end{figure}

The first example compares the influence of different input features(permeability matrices $K$ and transmissibility matrices $T$) on the prediction accuracy of the proxy model under the same network model. For the test case shown in \autoref{fig:perm-press}, there are four producing wells at four corners. We generate 2100 different permeability samples and compute the pressure fields at 150 time steps. There are 1500 samples in the training set and 600 samples in the test set. The pressure fields at the last time step with the highest variation are used for comparison.

\begin{figure}[H]
    \centering
    \subfloat[]{\includegraphics[width=0.5\linewidth]{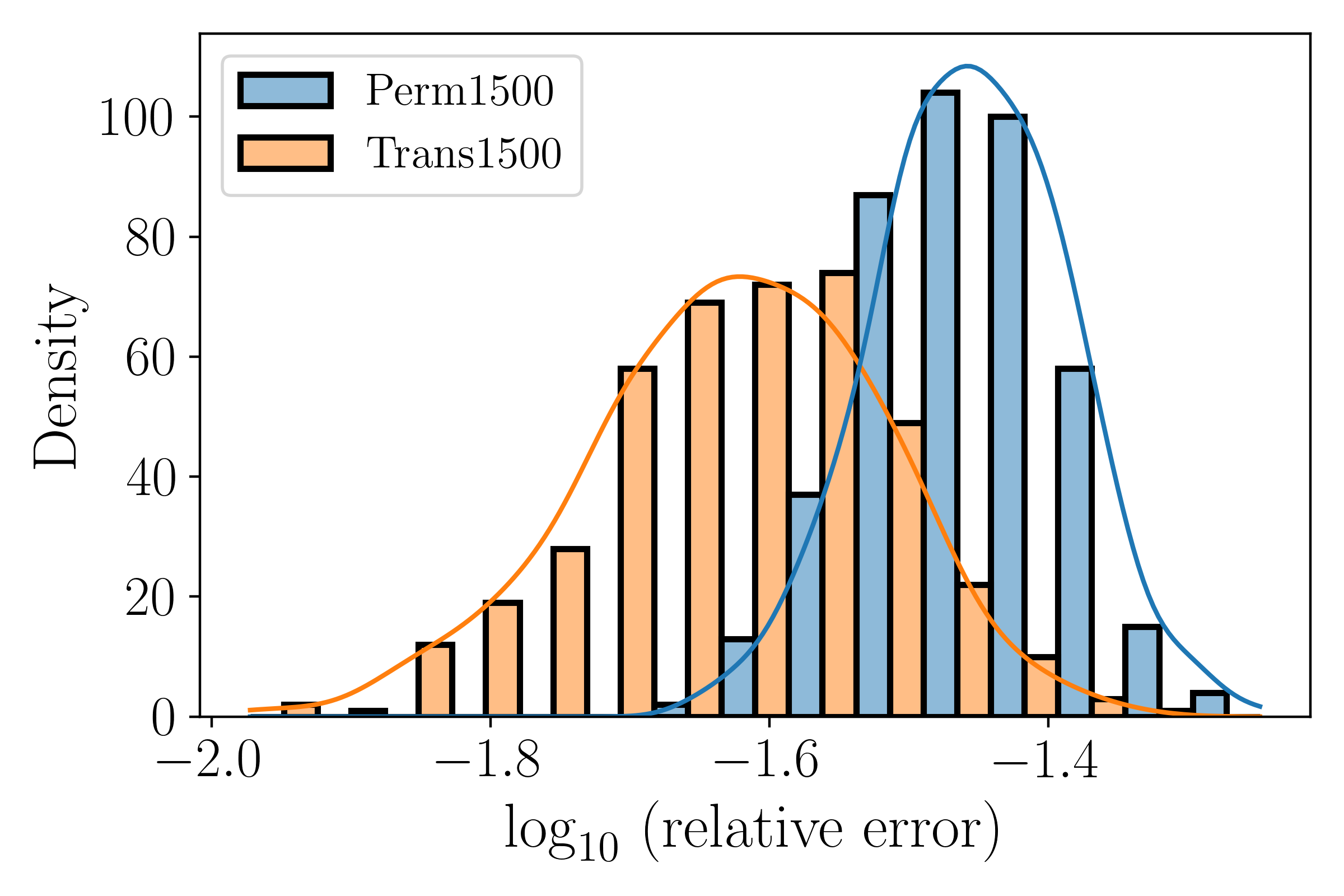}}
    \subfloat[]{\includegraphics[width=0.5\linewidth]{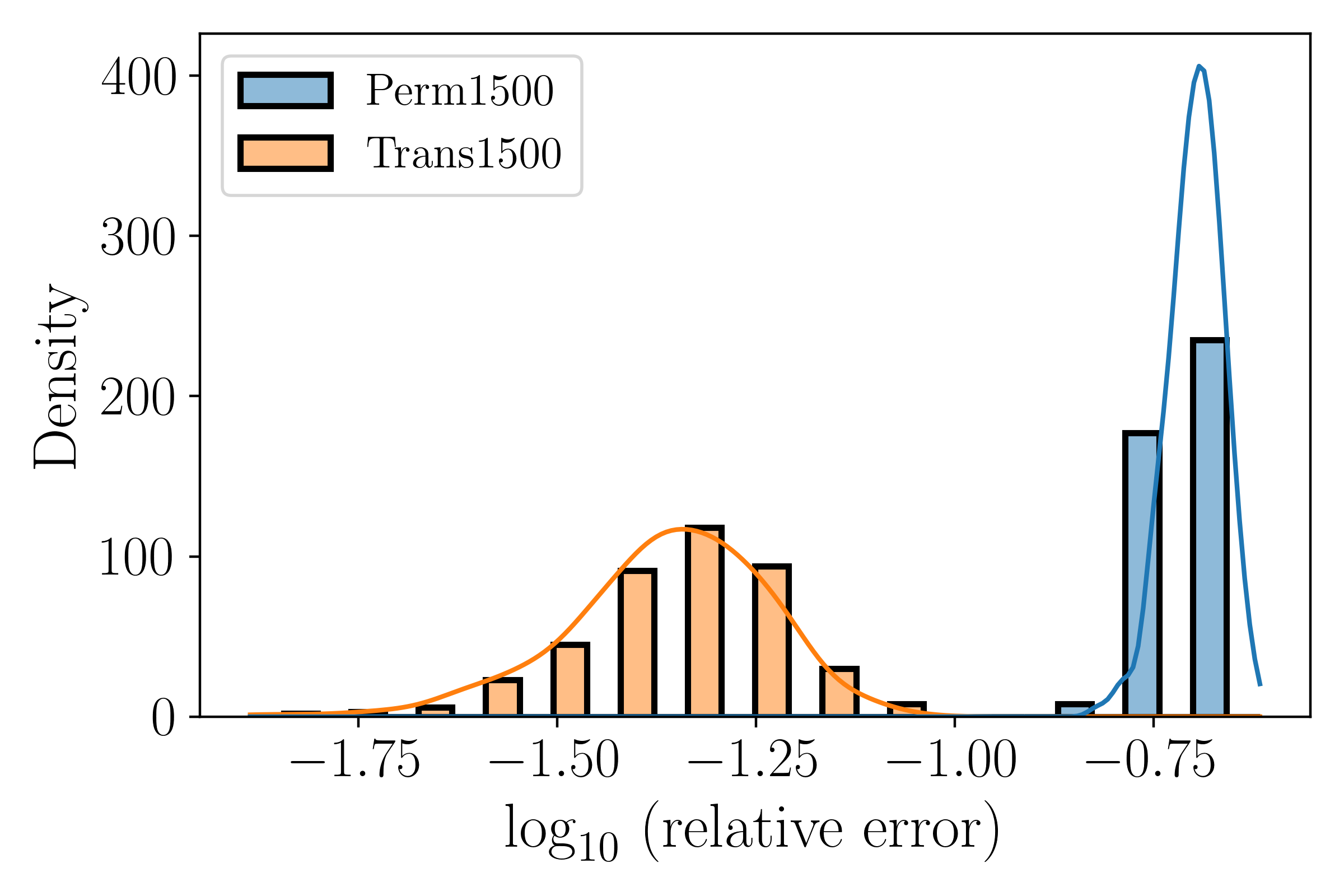}}
    \caption{The ARUnet surrogates are trained with the same samples but different input feature to predict the pressure field. (a) $\log_{10}{\text{relative error}}$ distribution on all grid cells in the test set. (b) $\log_{10}{\text{relative error}}$ on the well-block cell in the test set.}
    \label{fig:transVSperm}
\end{figure}

\begin{table}[htb]
    \centering
    \caption{Relative error using different input features}
    \begin{tabular}{l l l l}
    \toprule
      Input data & Relative error of producing well & Mean relative error (MRE) \\
    \midrule
      permeability   &  $2.01 \times10^{-1}$ & $3.54 \times10^{-2}$\\
      transmissibility & $4.53\times 10^{-2}$  & $2.44 \times10^{-2}$ \\
    \bottomrule
    \end{tabular}
    \label{tab:trans-perm}
\end{table}

As shown in \autoref{fig:transVSperm} and \autoref{tab:trans-perm}, the prediction accuracy is affected by inputs considerably. The neural network using transmissibility as inputs achieves better generalization ability compared to that using permeability. Specifically, the relative error at well-blocks (grid cells containing wells) decreases by 77\% to $4.53 \times 10^{-2}$, while the MRE across all grid cells shows a 31\% reduction to $2.44 \times10^{-2}$. 

\begin{figure}[H]
    \centering
    \subfloat[]{\includegraphics[width=0.33\textwidth]{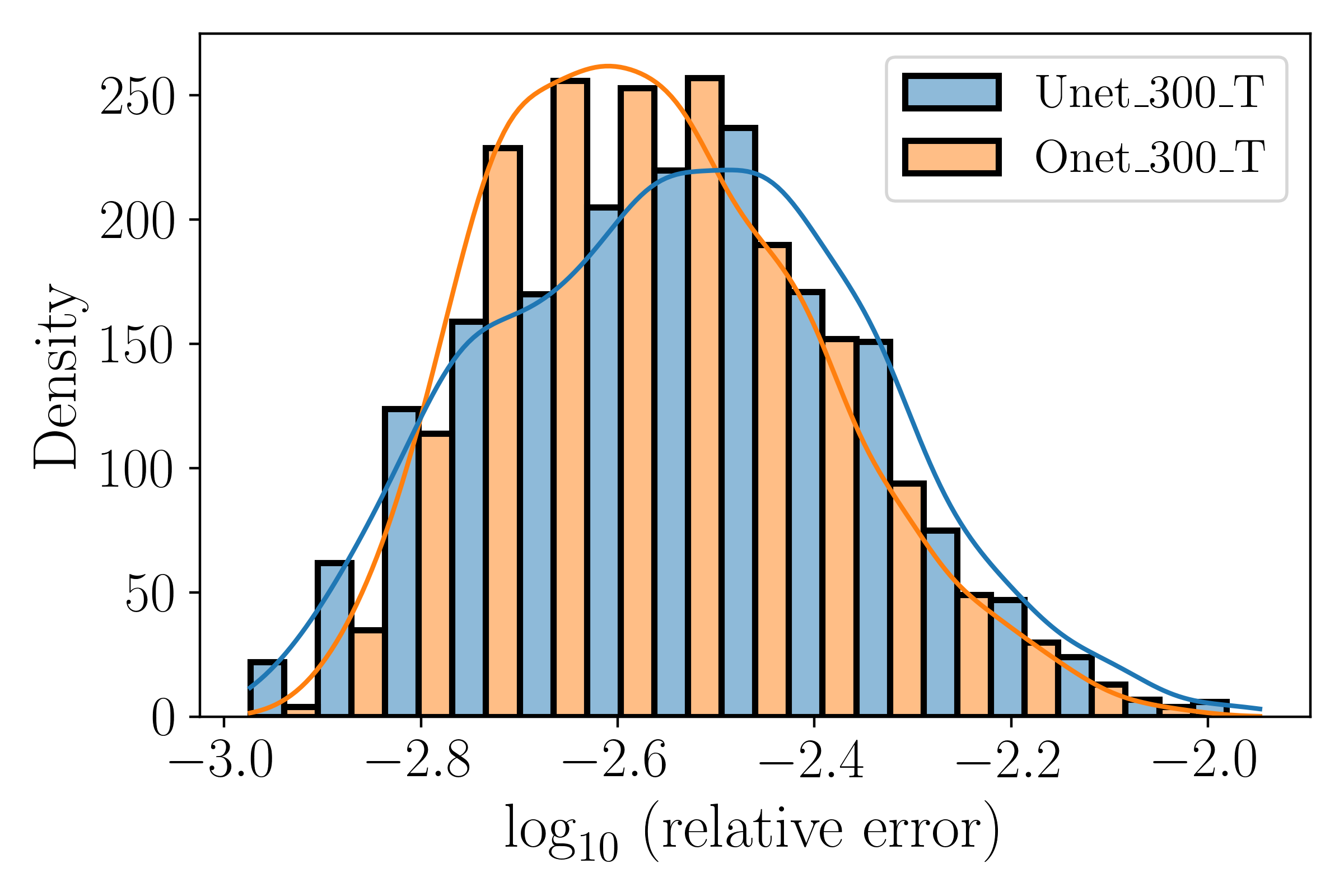}}
    \subfloat[]{\includegraphics[width=0.33\textwidth]{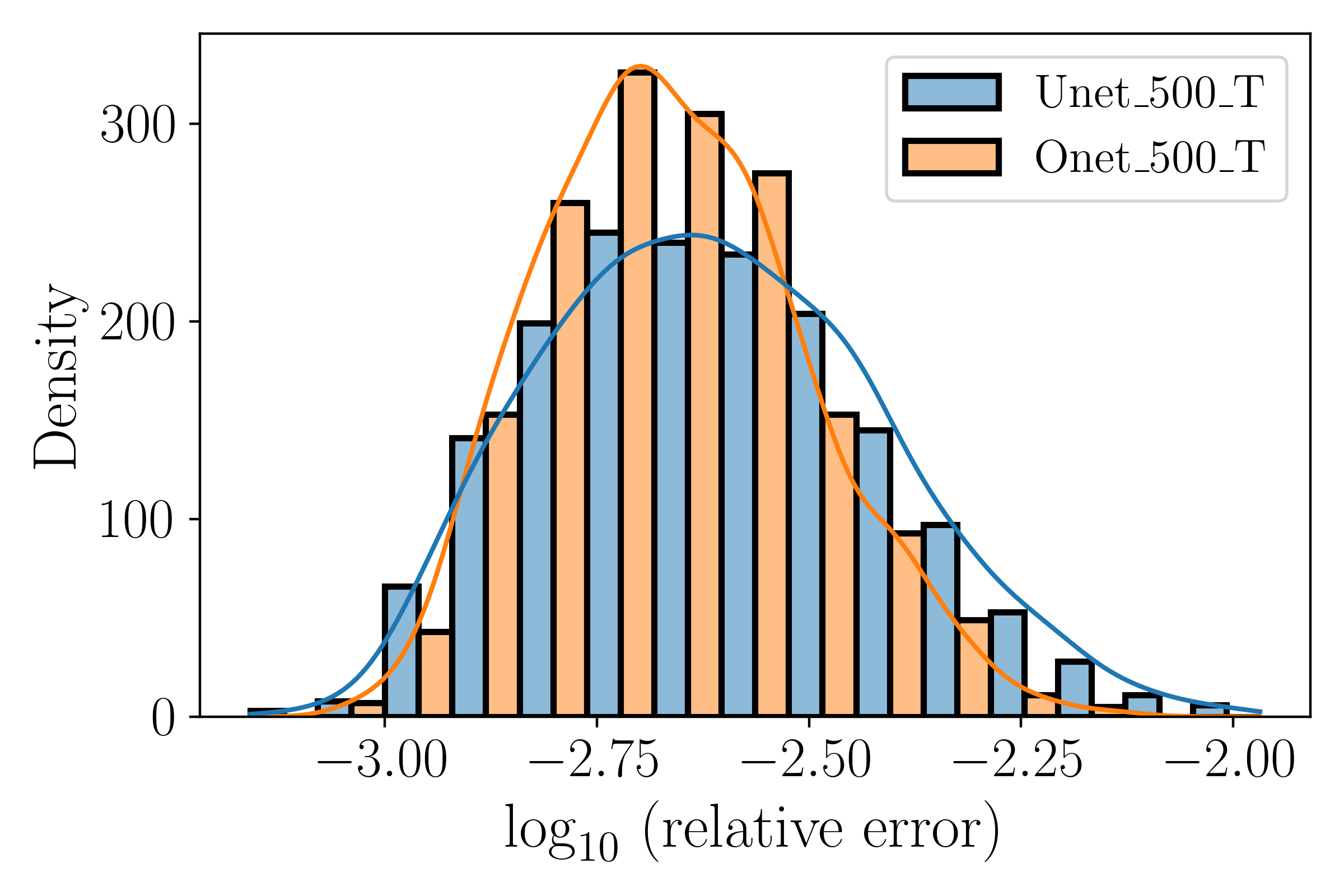}}
    \subfloat[]{\includegraphics[width=0.33\textwidth]{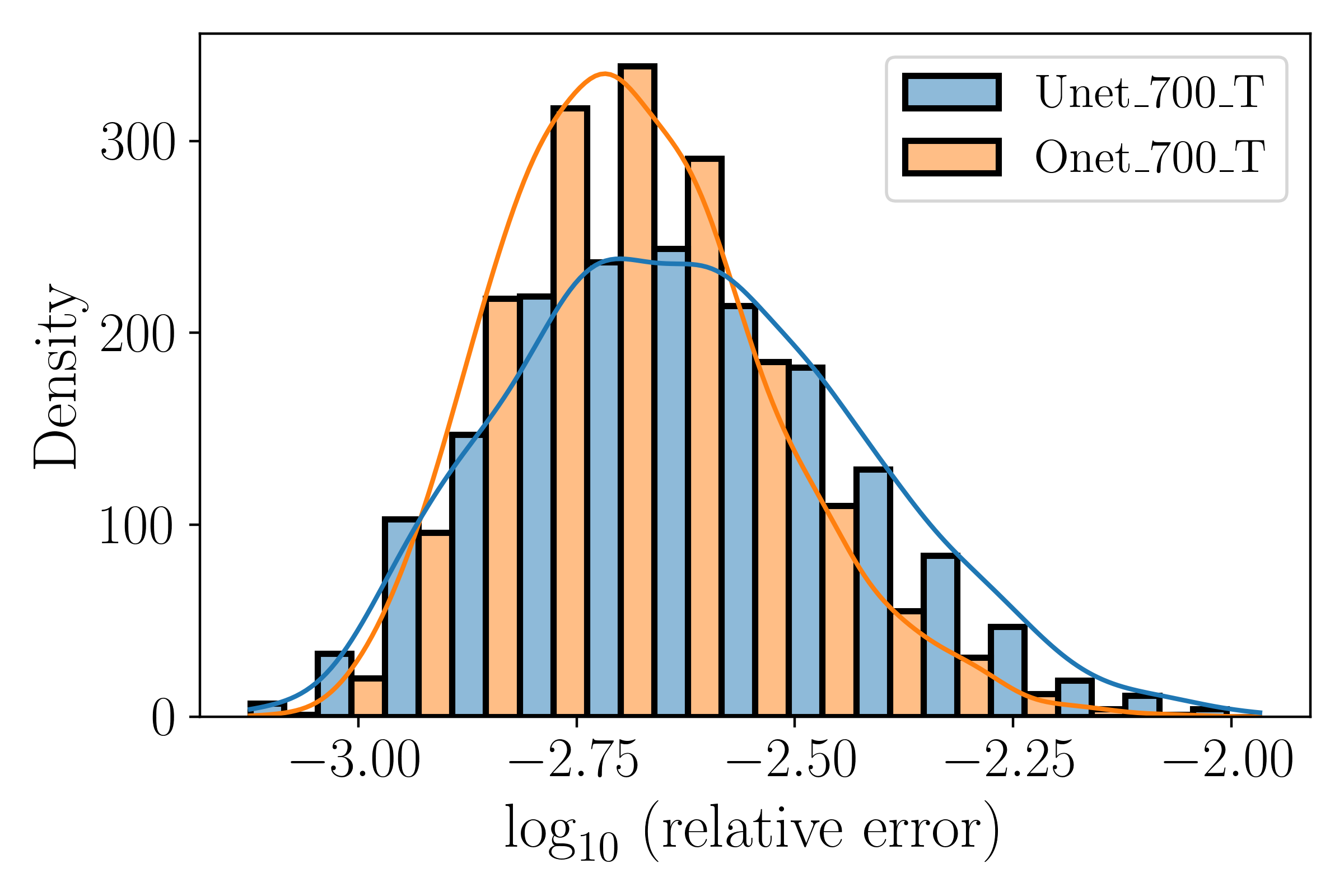}}

    \caption{$\log_{10}$ (MRE) distribution of pressure field prediction using ARUnet and AROnet with different samples. (a): 300 samples. (b): 500 samples. (c): 700 samples.}
    \label{fig:2d_net_sturcture_contrast}
\end{figure}

The second case compared the performance of ARUnet and AROnet, the training parameters are shown in \autoref{tab:nn-params}. We use the same training parameters and training data for both AROnet and ARUnet. we compared their performance in three different sample sizes: 300, 500, and 700 over 150 time steps. The results are shown on \autoref{tab:nn-params} and \autoref{fig:2d_net_sturcture_contrast}. 

\begin{table}[ht]
    \centering
    \caption{Configuration of Training Parameters for Two-Dimensional Pressure Field Prediction}
    \begin{tabular}{l l}
    \toprule
      Parameters   & Value \\
    \midrule
      Epochs   &  10 \\
      Iterations & 120 \\
      Optimizer  & Adam \\
      Sample interval & 10 \\
      Sample times & 10 \\
      Initial samples & 100 \\
      Sample number each time & 20, 40, 60 \\
      Learning rate & $5 \times 10^{-3}$ \\
    \bottomrule
    \end{tabular}
    \label{tab:nn-params}
\end{table}

\begin{table}[htb]
    \centering
    \caption{Results of the two-dimensional single-phase flow test case}
    \begin{tabular}{l l l l}
    \toprule
      samples  & NN structure  & mean relative error \\
    \midrule
      300   & ARUnet  & $3.22 \times 10^{-3} $ \\
      300   & AROnet  & $2.94 \times 10^{-3} $ \\
      500   & ARUnet & $2.77 \times 10^{-3} $ \\
      500   & AROnet  & $2.33 \times 10^{-3} $ \\
      700   & ARUnet  & $2.66 \times 10^{-3} $ \\
      700   & AROnet  & $2.21 \times 10^{-3} $ \\
    \bottomrule
    \end{tabular} 
    \label{tab:2dim-result}
\end{table}

\begin{figure}[ht]
    \centering
    \subfloat[]{\includegraphics[width=0.45\textwidth]{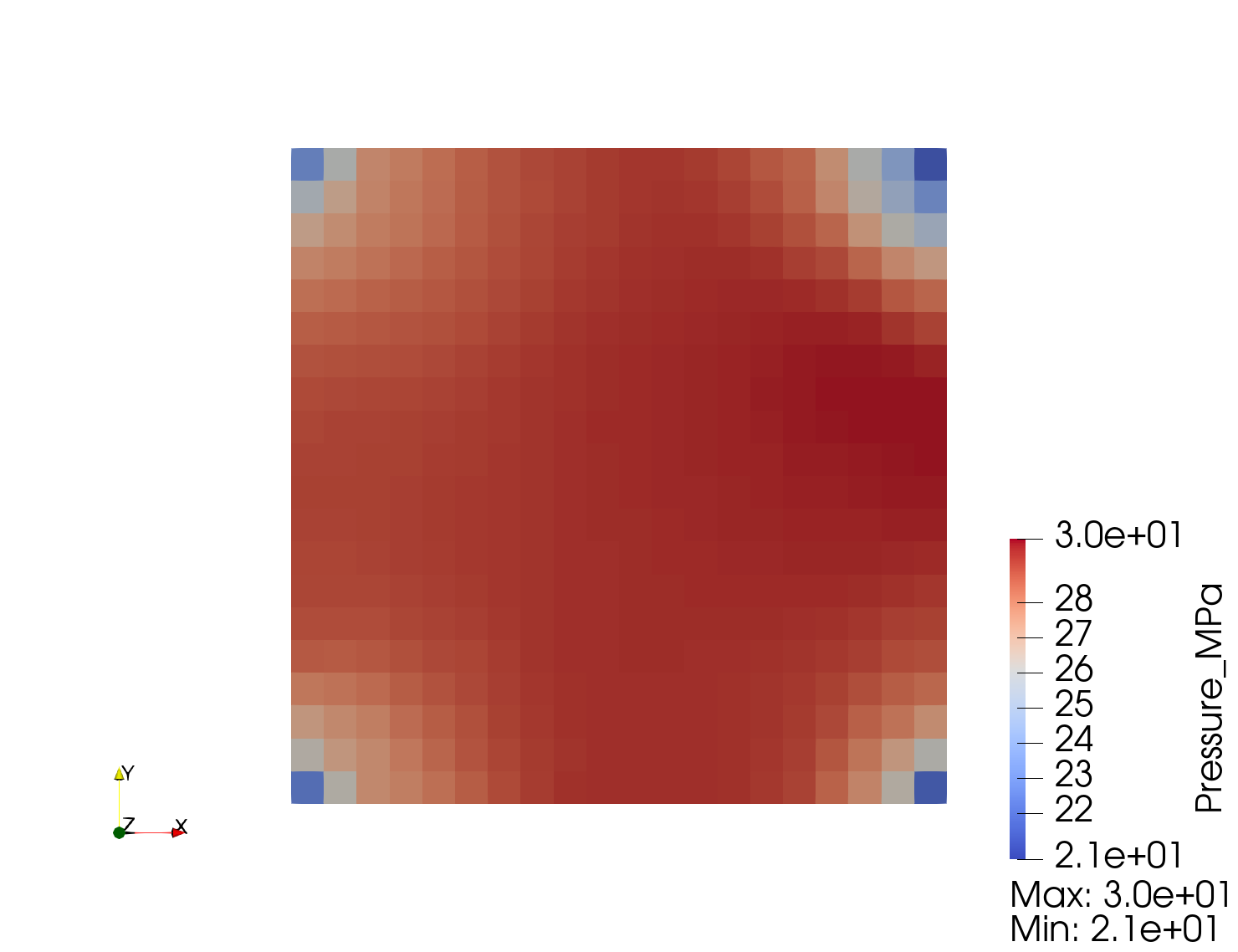}}
    \subfloat[]{\includegraphics[width=0.45\textwidth]{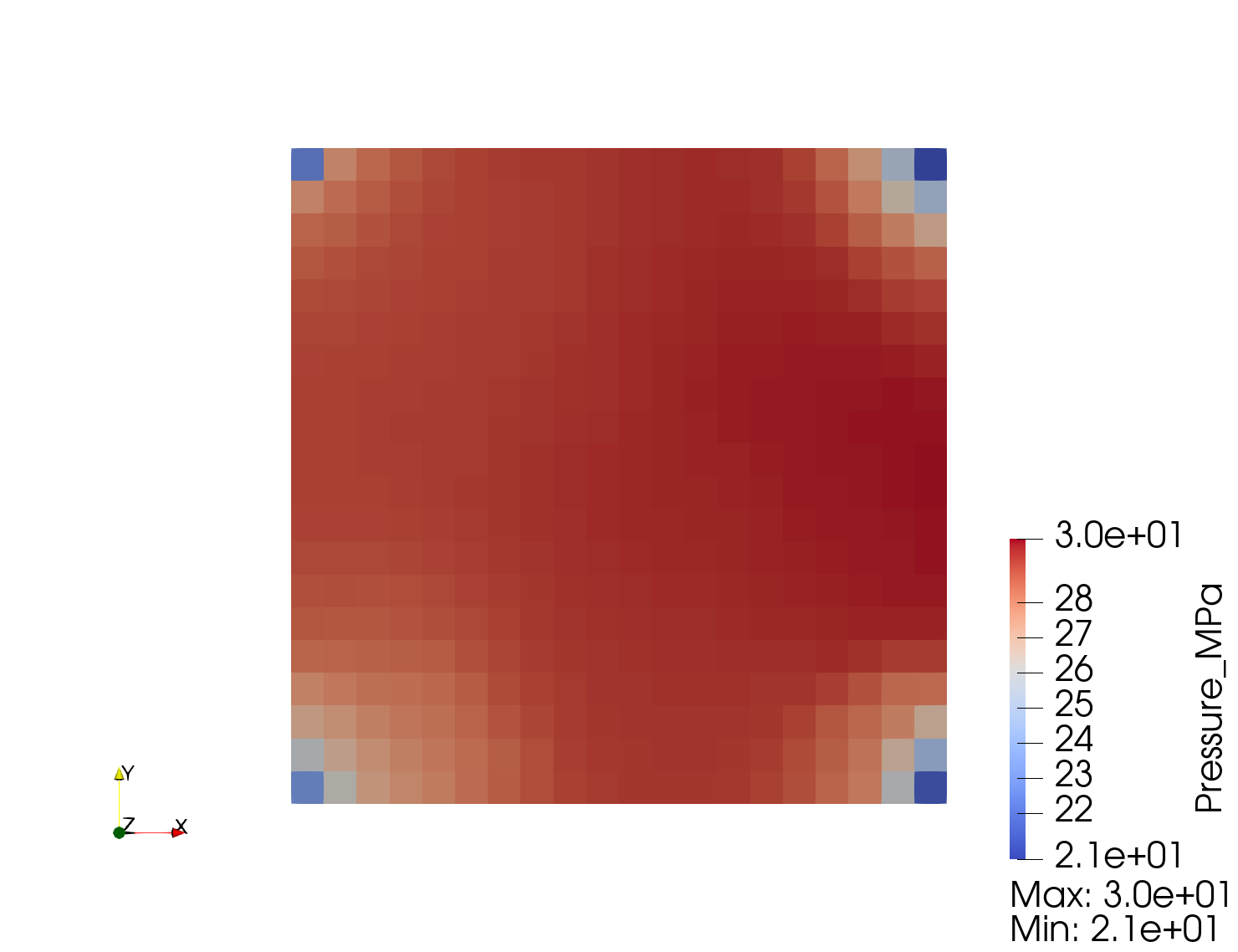}}\\

    \subfloat[]{\includegraphics[width=0.45\textwidth]{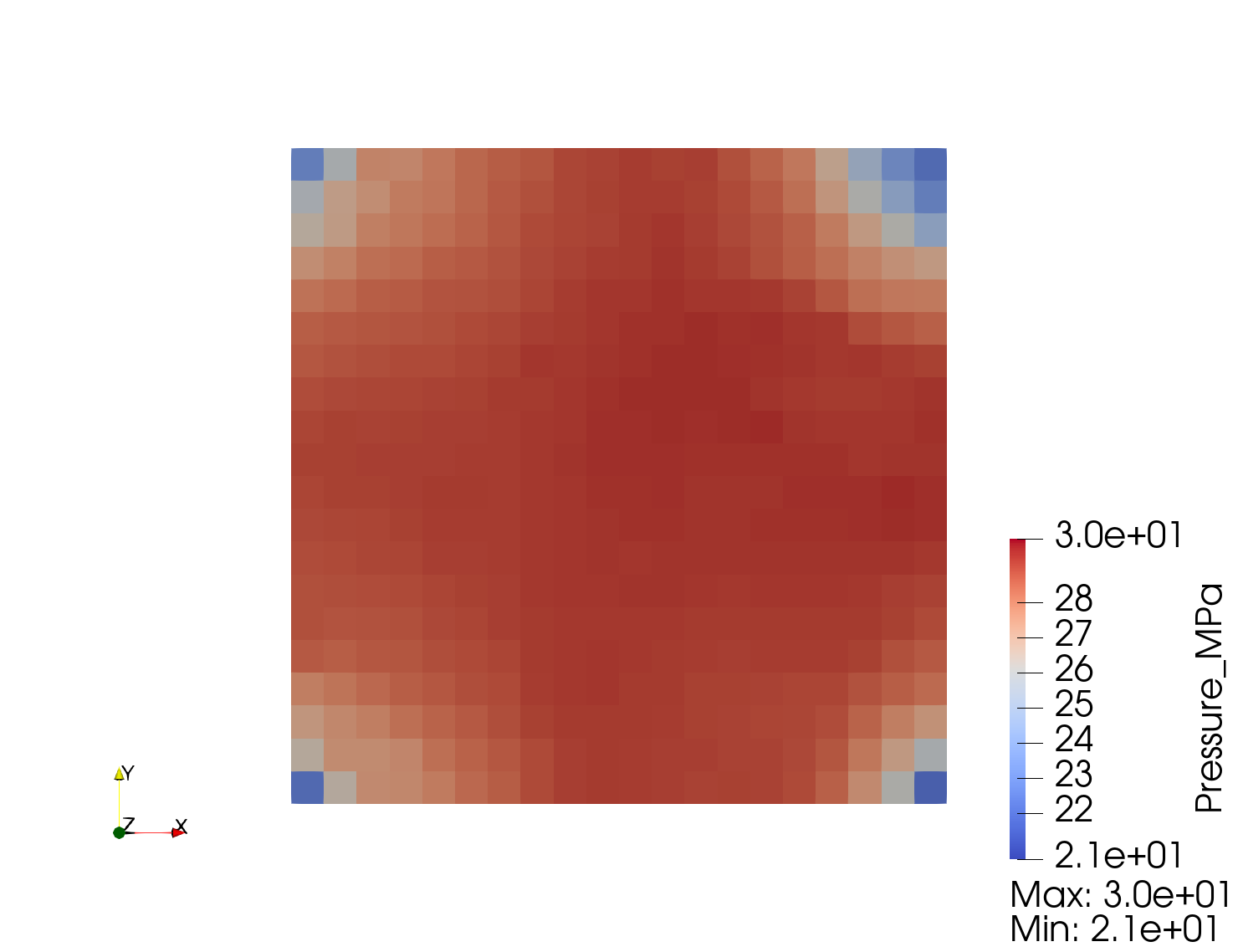}}
    \subfloat[]{\includegraphics[width=0.45\textwidth]{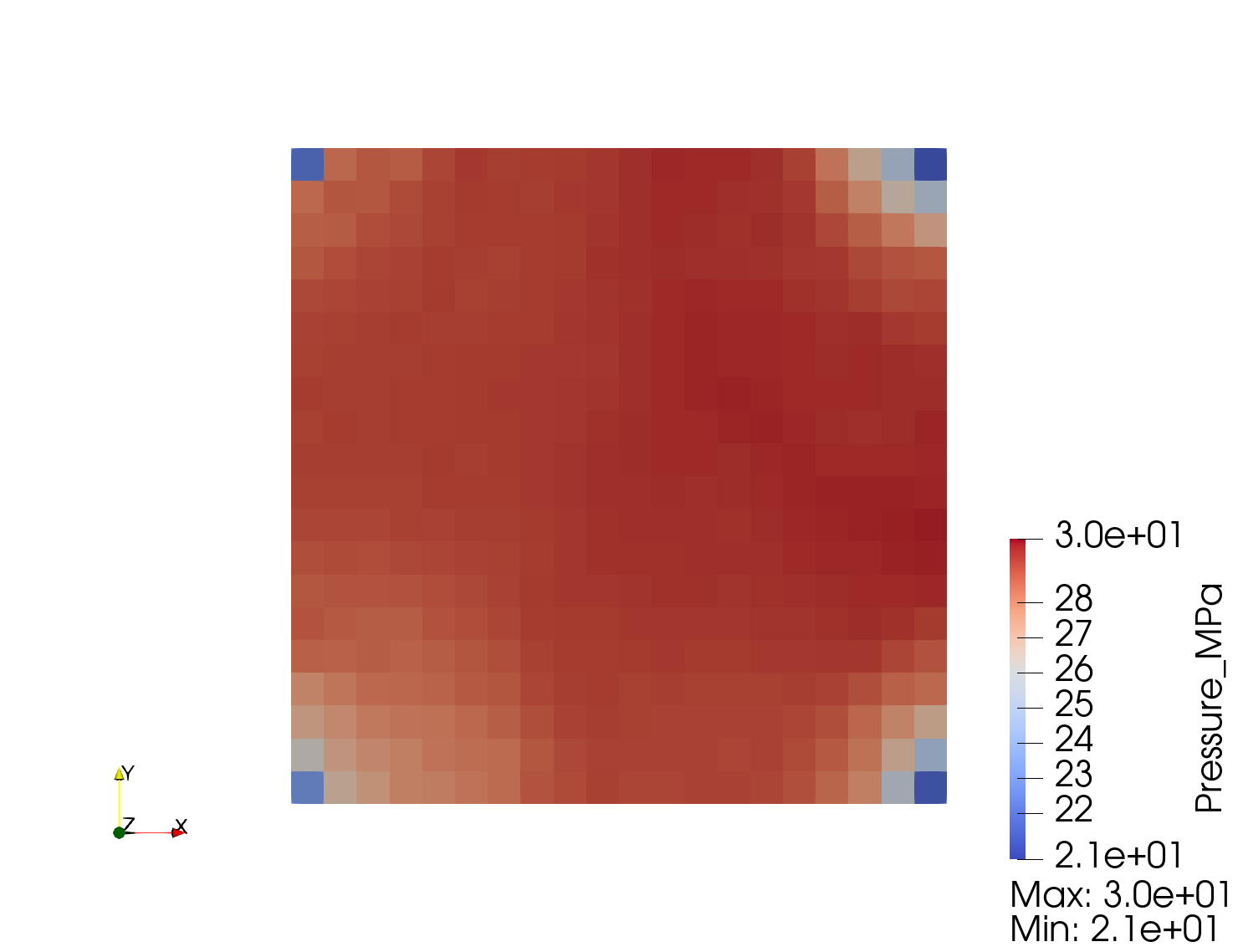}}
    \caption{(a), (b): Ground truth of two different pressure fields of test samples simulated using FVM . (c), (d): Corresponding pressure fields predicted by AROnet}
    \label{fig:2d-press}
\end{figure}

As demonstrated in \autoref{fig:2d_net_sturcture_contrast} and \autoref{tab:2dim-result}, AROnet consistently outperforms ARUnet in pressure field prediction accuracy. This validates the advantage of AROnet's branch-trunk architecture, which explicitly encodes spatial correlations in heterogeneous parameter fields through multiplicative feature interactions. This advantage is also demonstrated in \autoref{sec:2phase-val}.
In \autoref{fig:2d-press}, two samples from the test set are visualized to validate the predicted pressure field by AROnet against the ground truth simulated using FVM.

\subsection{3D Single-Phase Darcy Flow}\label{sec:3D-val}

\begin{figure}[H]
    \centering
    \subfloat[]{\includegraphics[width=0.45\linewidth]{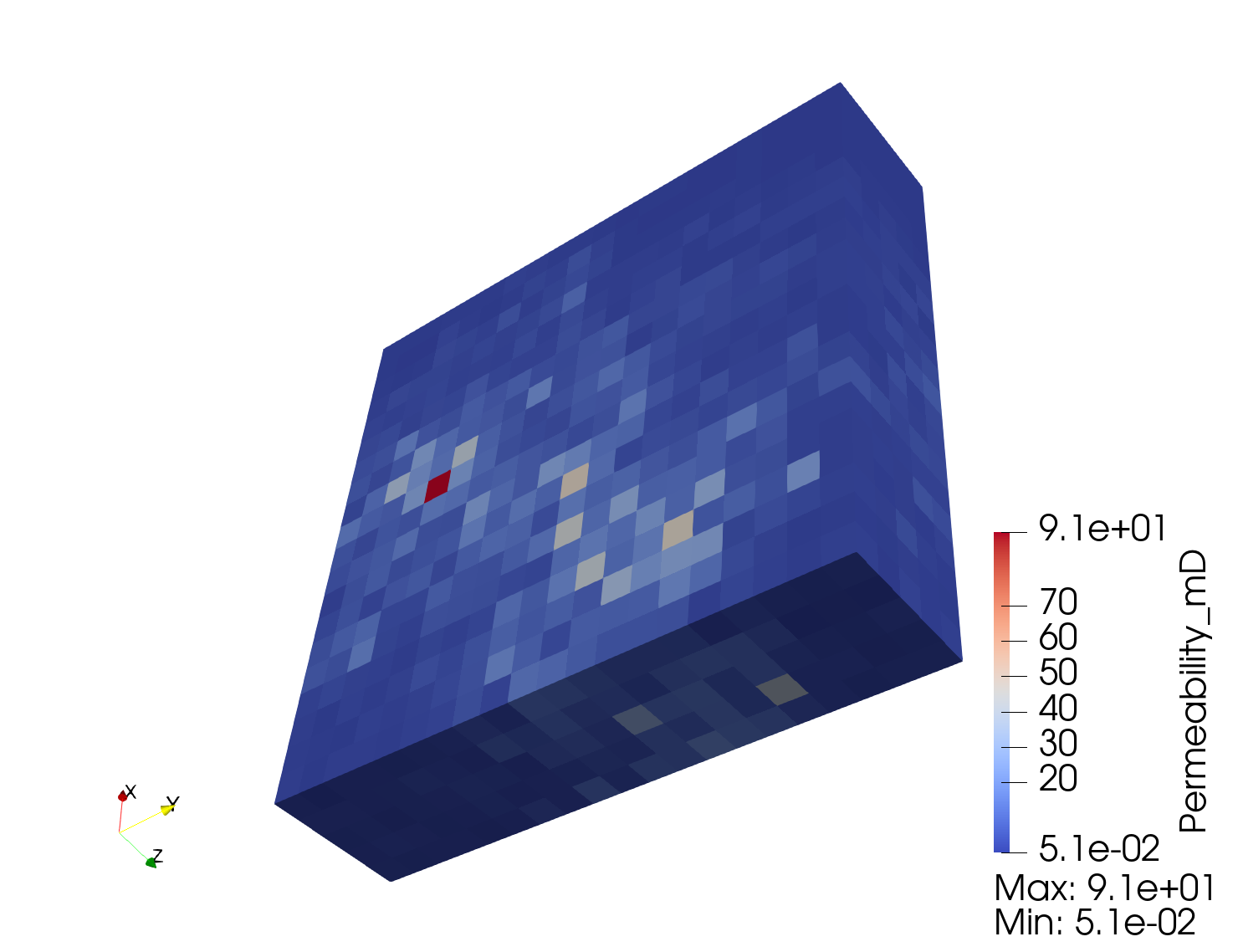}}
    \subfloat[]{\includegraphics[width=0.45\linewidth]{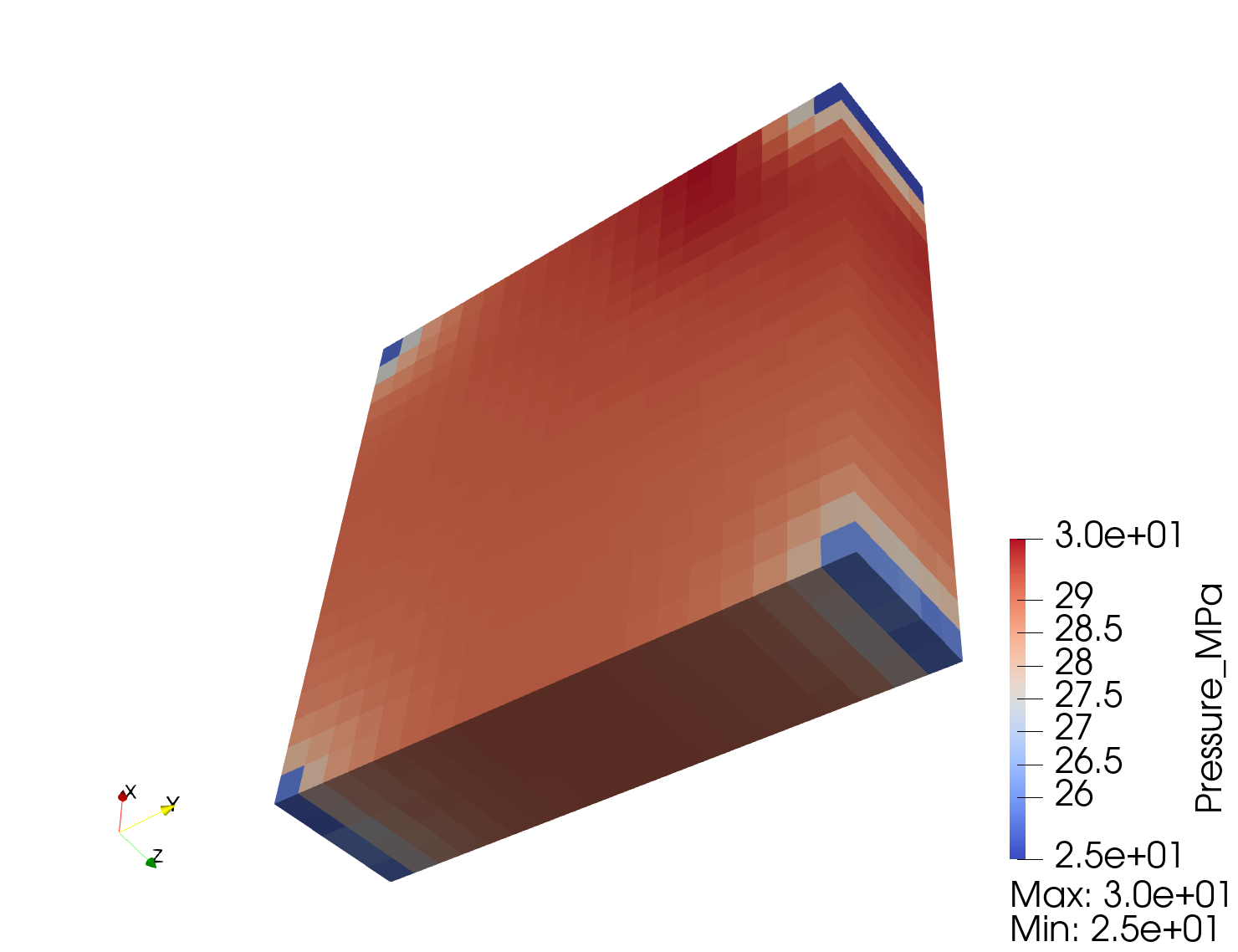}}
    
    \caption{(a): Heterogeneous 3D permeability field generated by Gaussian sequential simulation. (b): Corresponding pressure field computed by FVM at the last time step.}
    \label{fig:3d-perm-press}
\end{figure}

To accommodate 3D flow fields (see \autoref{fig:3d-perm-press}), the model architecture was extended by reshaping input tensors from [batch size, 4, H, W] to [batch size, 6 $\times$ d, H, W], where d denotes the depth dimension of the input parameter field. GMM is implemented for adaptive sampling to improve prediction accuracy shown in \autoref{tab:3d-result} while maintaining computational efficiency. Key implementation details include
\begin{itemize}
    \item Network architecture and hyperparameters (see \autoref{tab:nn_3d-params}).
    \item Depth-aware feature concatenation for 3D fields.
    \item GMM-driven adaptive sampling.
\end{itemize}

\begin{figure}[H]
	\centering
	\subfloat[]{\includegraphics[width=0.5\linewidth]{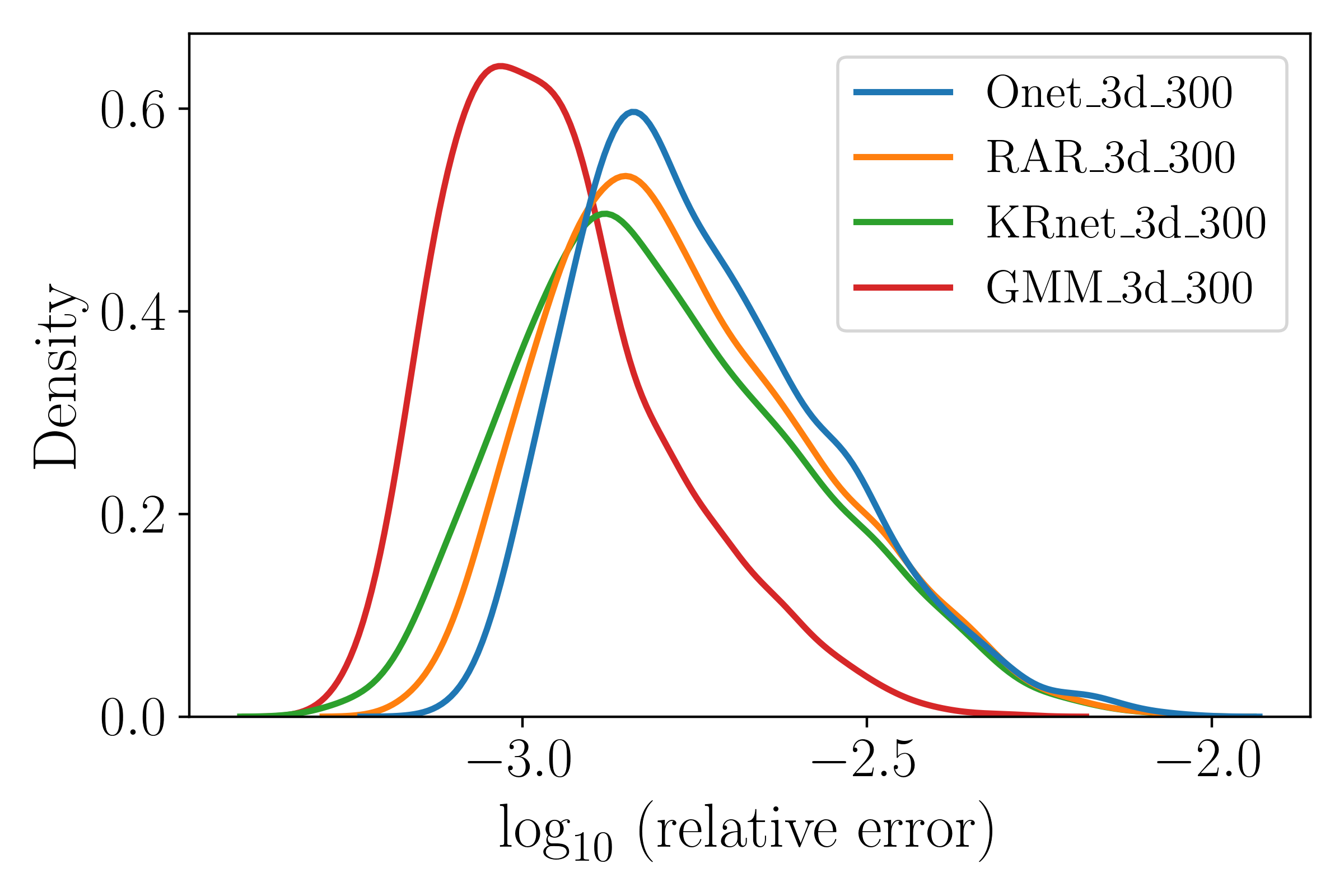}}
	\subfloat[]{\includegraphics[width=0.5\linewidth]{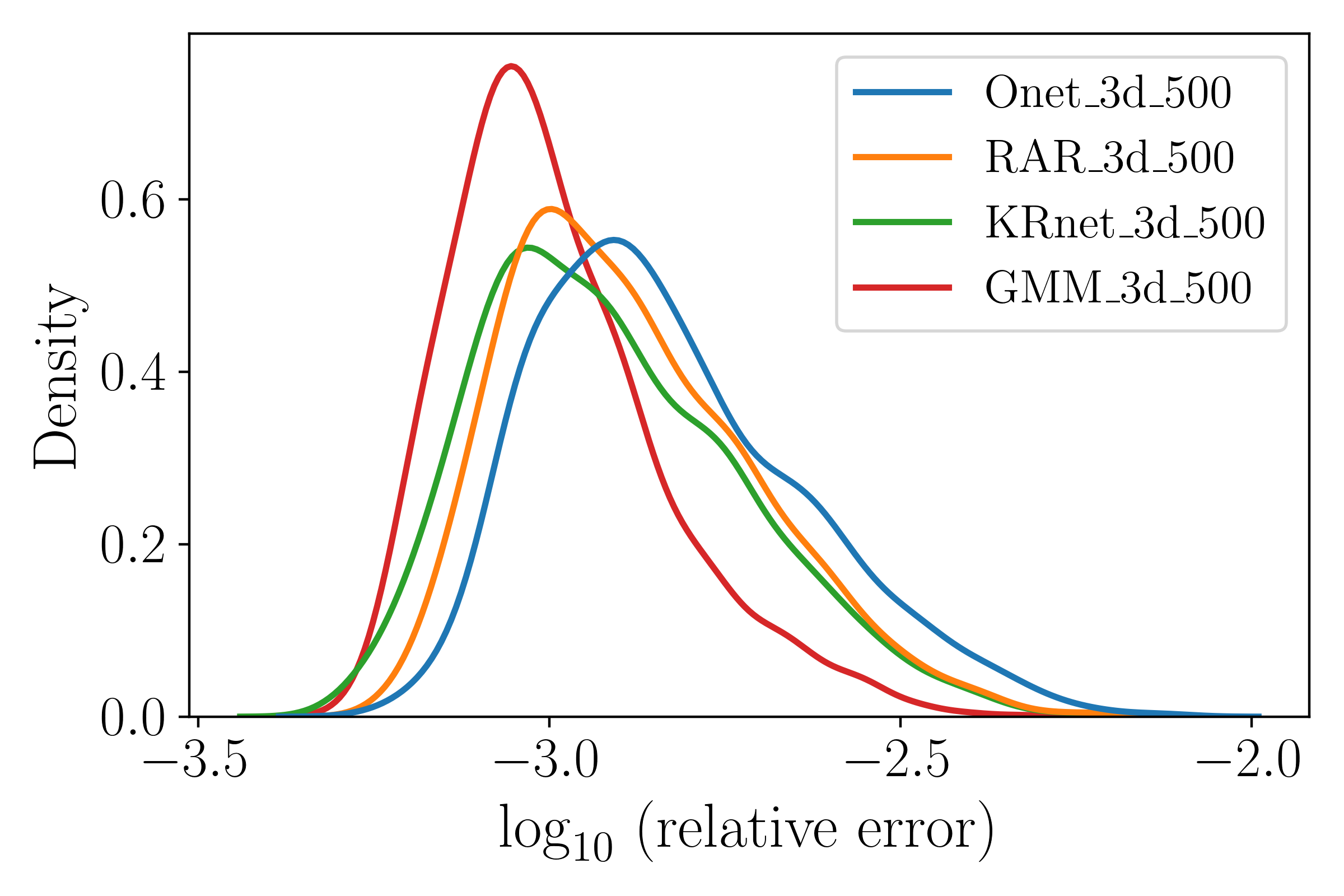}}
	\caption{ (a) Mean relative error distribution in test set corresponding to a training set with 300 samples, where curve ResAS300 is for RAR, KRnet300 is for KRnet, Unet300 is for random sampling and GMM300 is for GMM. (b) Results corresponding to a training set with 500 samples.}
	\label{fig:AS-effect}
\end{figure}

Three adaptive sampling strategies are implemented for validation using the test case in \autoref{fig:3d-perm-press} with training set of size 300. The baseline is the conventional random sampling method. The comparison results are presented in \autoref{tab:AS-time}.
\autoref{fig:AS-effect} presents the relative error distributions for prediction of adaptive sampling methods. All models have the same 100 initial training samples. In the training process, new samples are generated and added to the training set. The results show that the GMM-based adaptive sampling method outperforms random sampling in terms of prediction accuracy, as evidenced by the lower MRE and MSE values. The RAR method also shows improvement over random sampling, but it is less effective than the GMM-based approach. The KRnet method, while computationally expensive, does not significantly outperform the other methods in terms of prediction accuracy.

\begin{table}[ht]
	\centering
	\caption{Average time consumption, MRE and MSE for prediction of different sampling methods with 300 samples in total for training.}
	\begin{tabular}{l l l l}
		\toprule
		Type   & sampling time & MRE & MSE \\
		\midrule
            Random & None & $2.03 \times10^{-3}$ & $8.61 \times10^{9}$\\ 
		RAR &  0.31s &  $1.93 \times10^{-3}$ & $8.28 \times10^{9}$\\
		KRnet & 625s &  $1.85 \times10^{-3}$ & $7.59 \times10^{9}$\\
            GMM & 0.60s  & $1.21 \times10^{-3}$ & $3.94 \times10^{9}$\\
		\bottomrule
	\end{tabular}
	
	\label{tab:AS-time}
\end{table}

\begin{table}[ht]
    \centering
    \caption{Settings for Training Parameters in Three Dimensional Pressure Field Prediction Tasks}
    \begin{tabular}{l l}
    \toprule
      Parameters   & Value \\
    \midrule
      Epochs   &  5 \\
      Iterations & 550 \\
      Optimizer  & Adam \\
      Sample interval & 50 \\
      Sample times & 10 \\
      Initial samples & 100 \\
      Sample number each time & 20, 40, 60 \\
      Learning rate & $5 \times 10^{-3}$ \\
    \bottomrule
    \end{tabular}
    \label{tab:nn_3d-params}
    \end{table}

    \begin{table}[H]
        \centering
        \caption{Results of three-dimensions single phase flow test case}
        \begin{tabular}{l l l l}
        \toprule
          samples  &   well relative error & mean relative error & adaptive sampling \\ 
        \midrule
          300   &  $6.12 \times 10^{-3}$ & $2.03 \times 10^{-3} $ & False \\
          300   &  $4.62 \times 10^{-3}$ & $1.21 \times 10^{-3} $ & True \\
          
          500   &  $3.65 \times 10^{-3}$ & $1.56 \times 10^{-3} $ & False \\
          500   &  $3.39 \times 10^{-3}$ & $9.65 \times 10^{-4} $ & True \\
          
          700   &  $3.41 \times 10^{-3}$ & $1.35 \times 10^{-3} $ & False \\
          700   &  $3.02 \times 10^{-3}$ & $8.30 \times 10^{-4} $ & True \\
        \bottomrule
        \end{tabular}
        \label{tab:3d-result}
    \end{table}

\begin{figure}[H]
    \centering
    \subfloat[]{\includegraphics[width=0.33\linewidth]{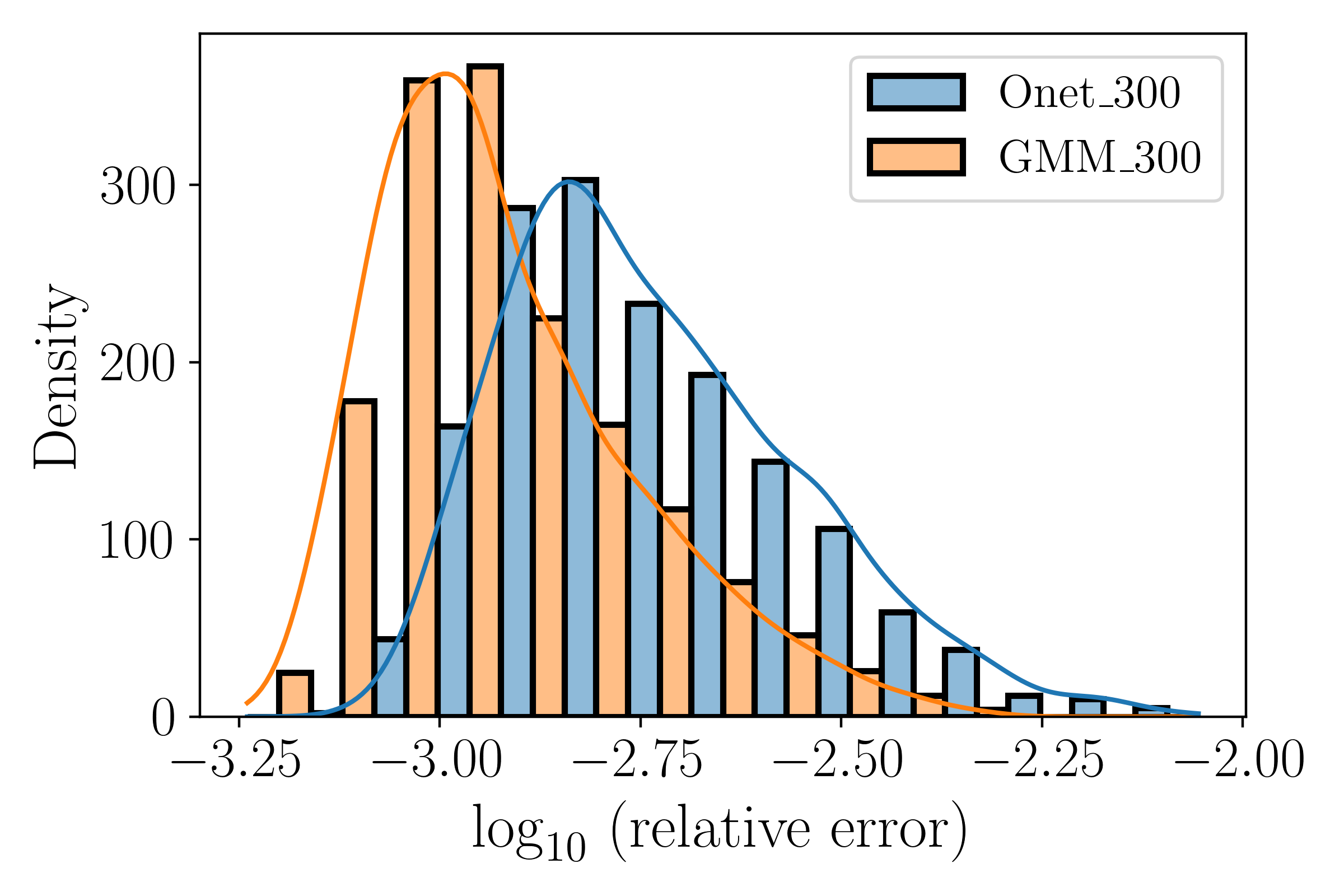}}
    \subfloat[]{\includegraphics[width=0.33\linewidth]{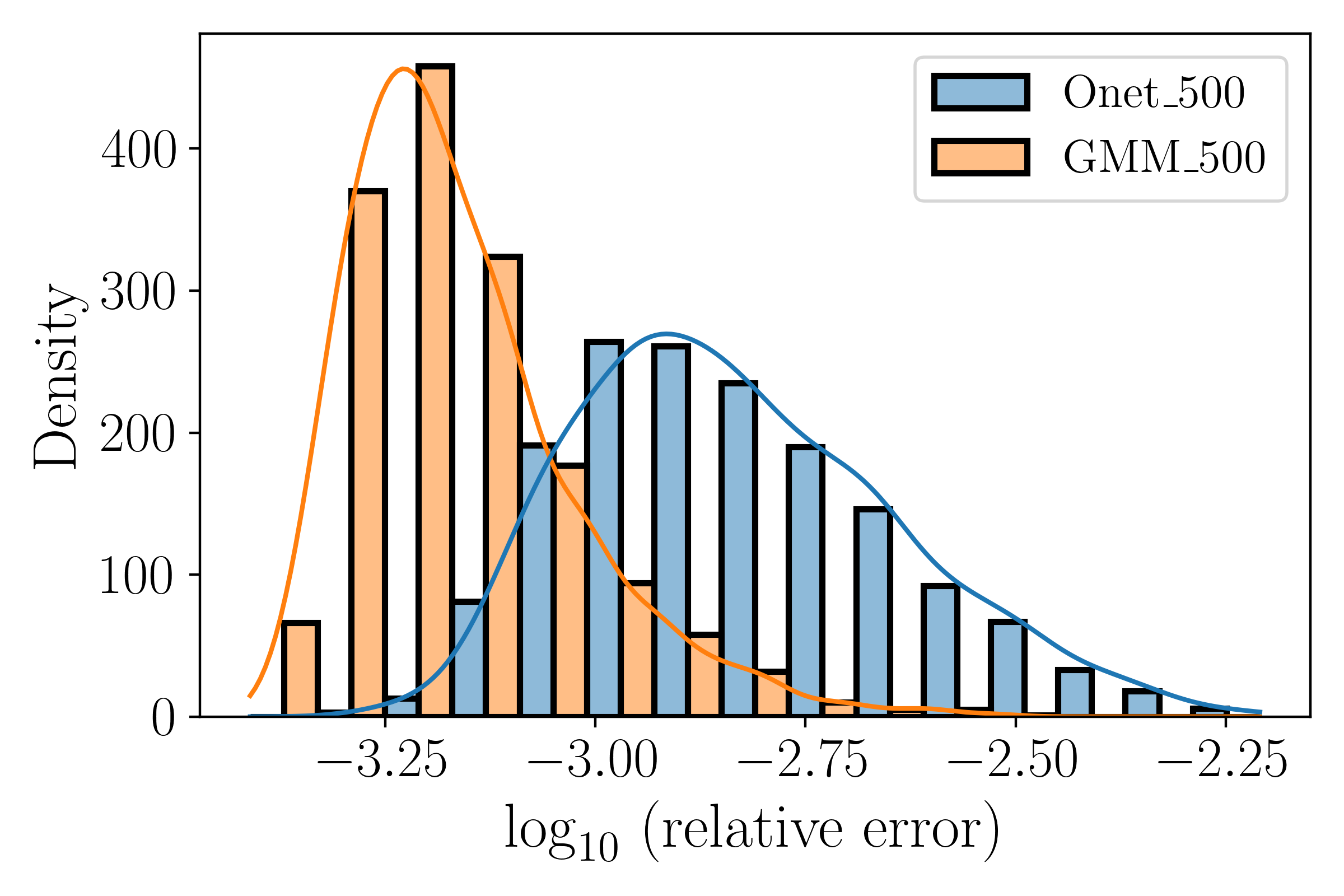}}
    \subfloat[]{\includegraphics[width=0.33\linewidth]{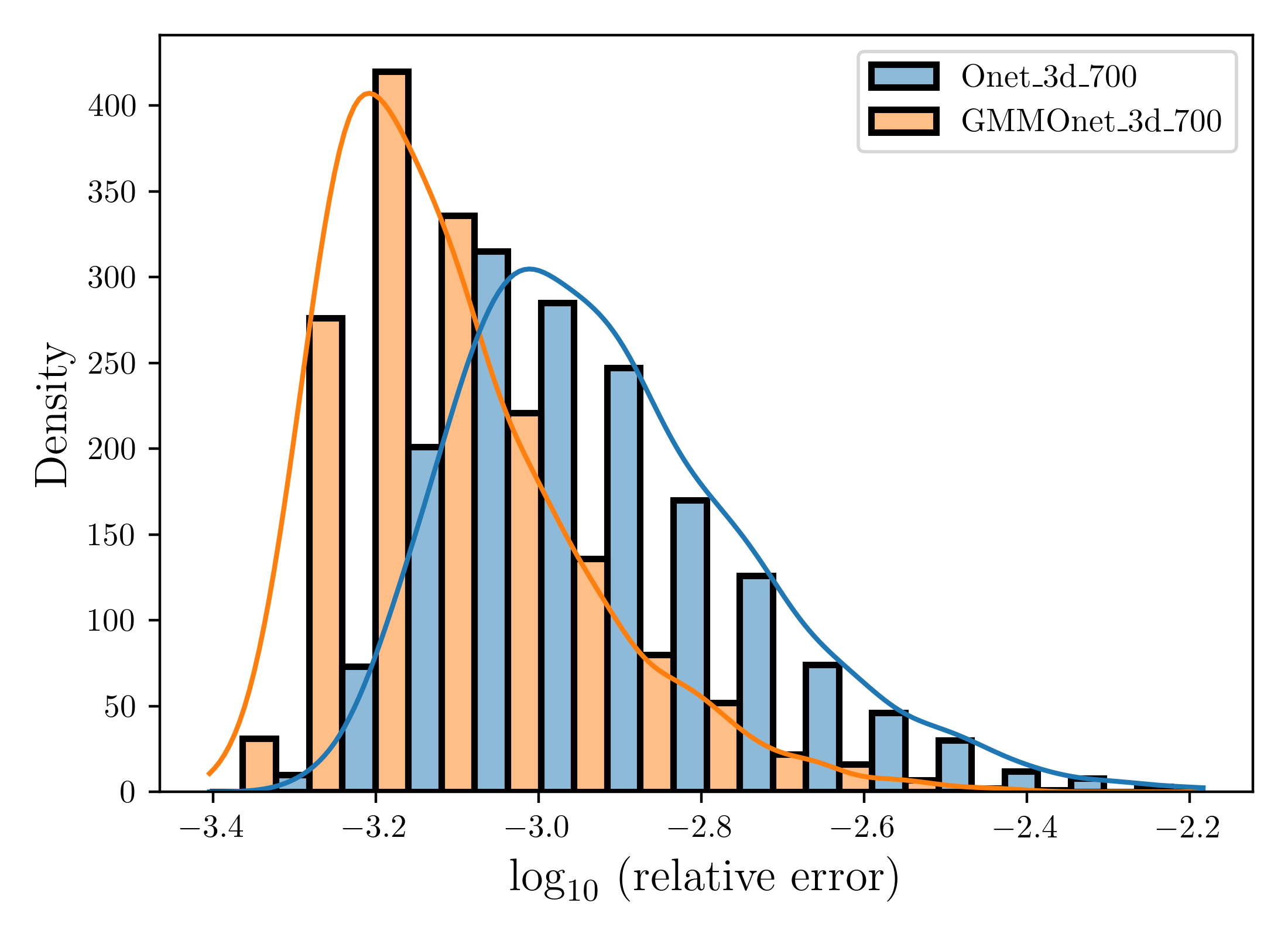}}
    \\
    \subfloat[]{\includegraphics[width=0.33\linewidth]{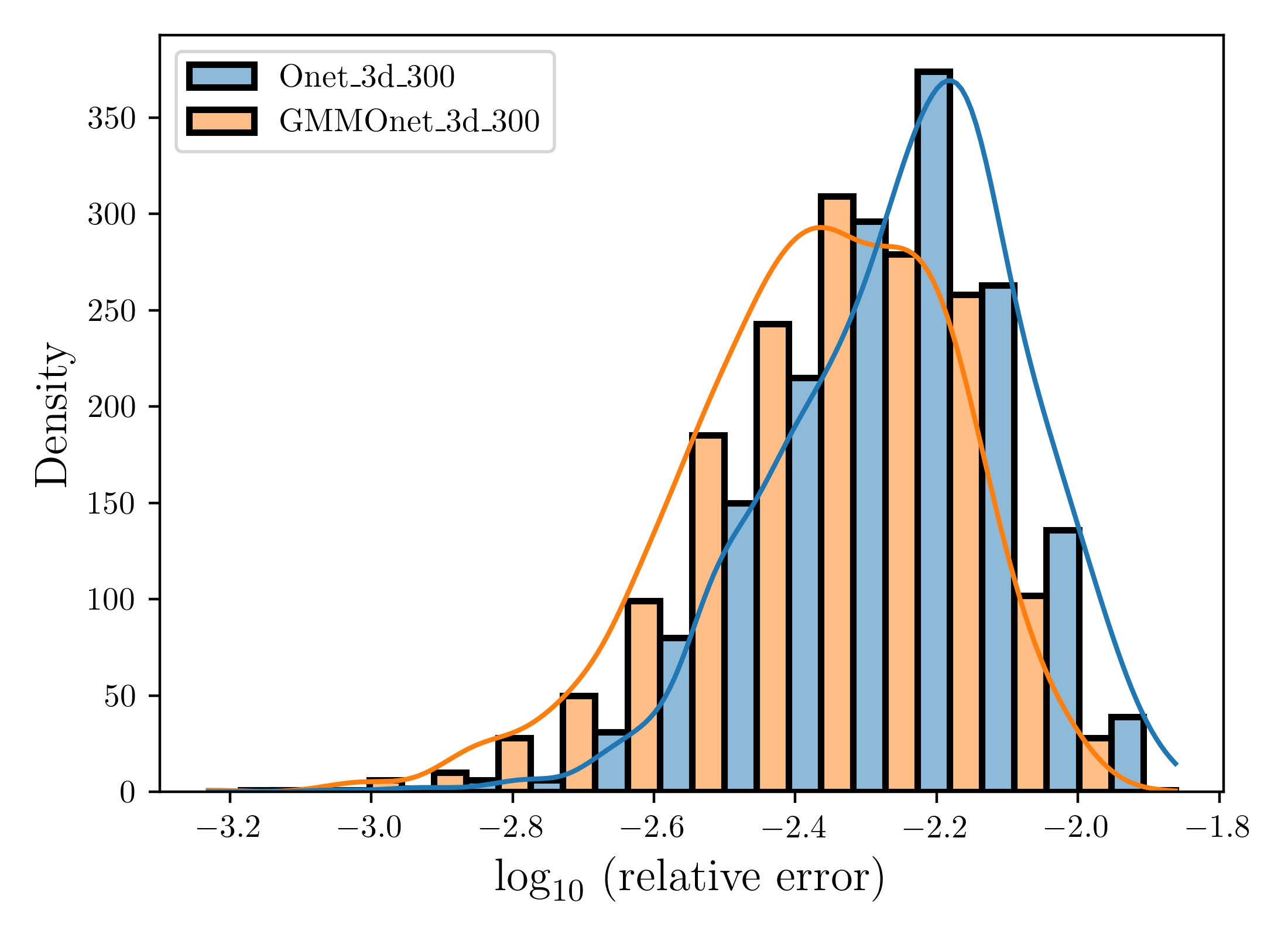}}
    \subfloat[]{\includegraphics[width=0.33\linewidth]{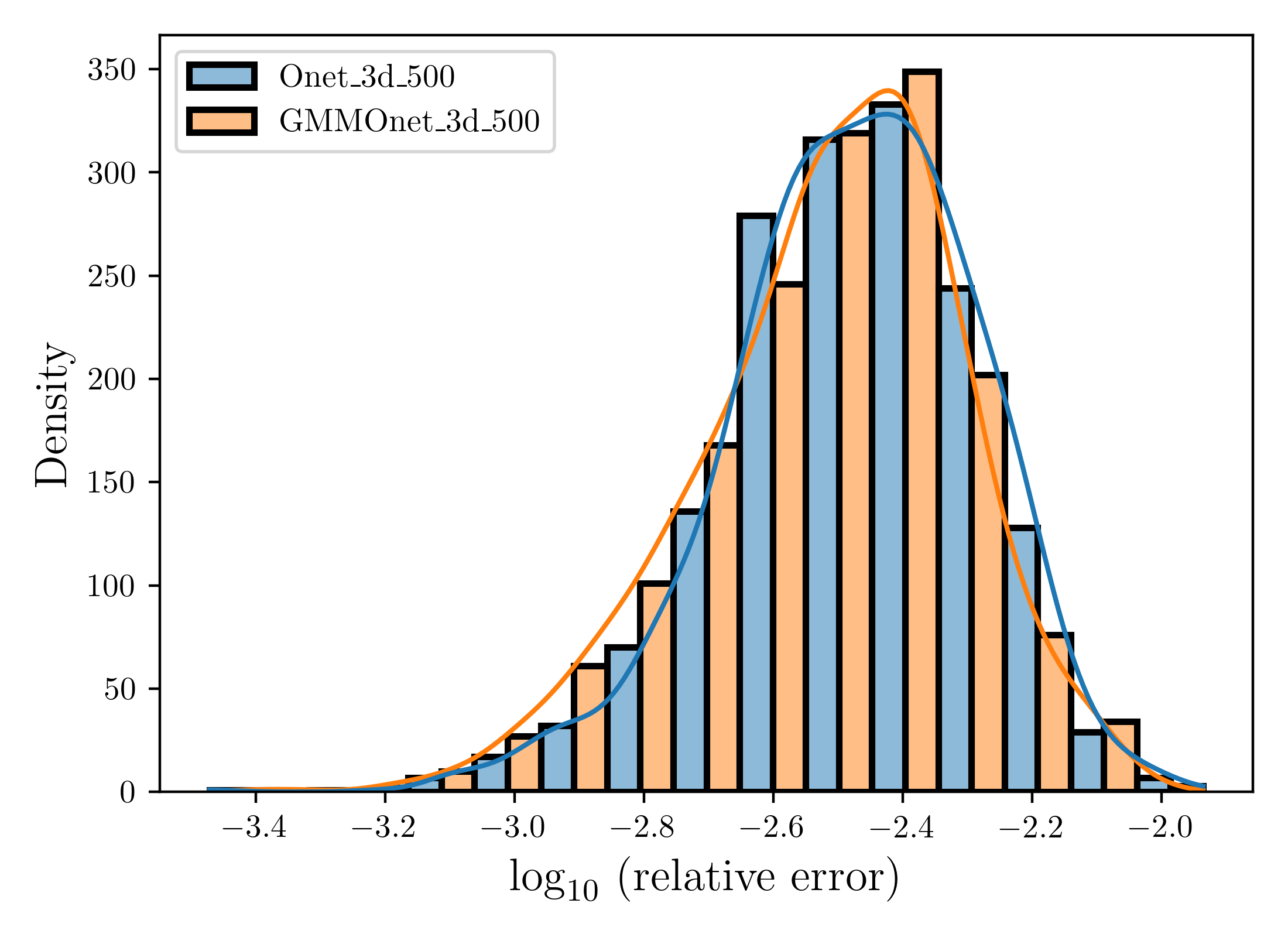}}
    \subfloat[]{\includegraphics[width=0.33\linewidth]{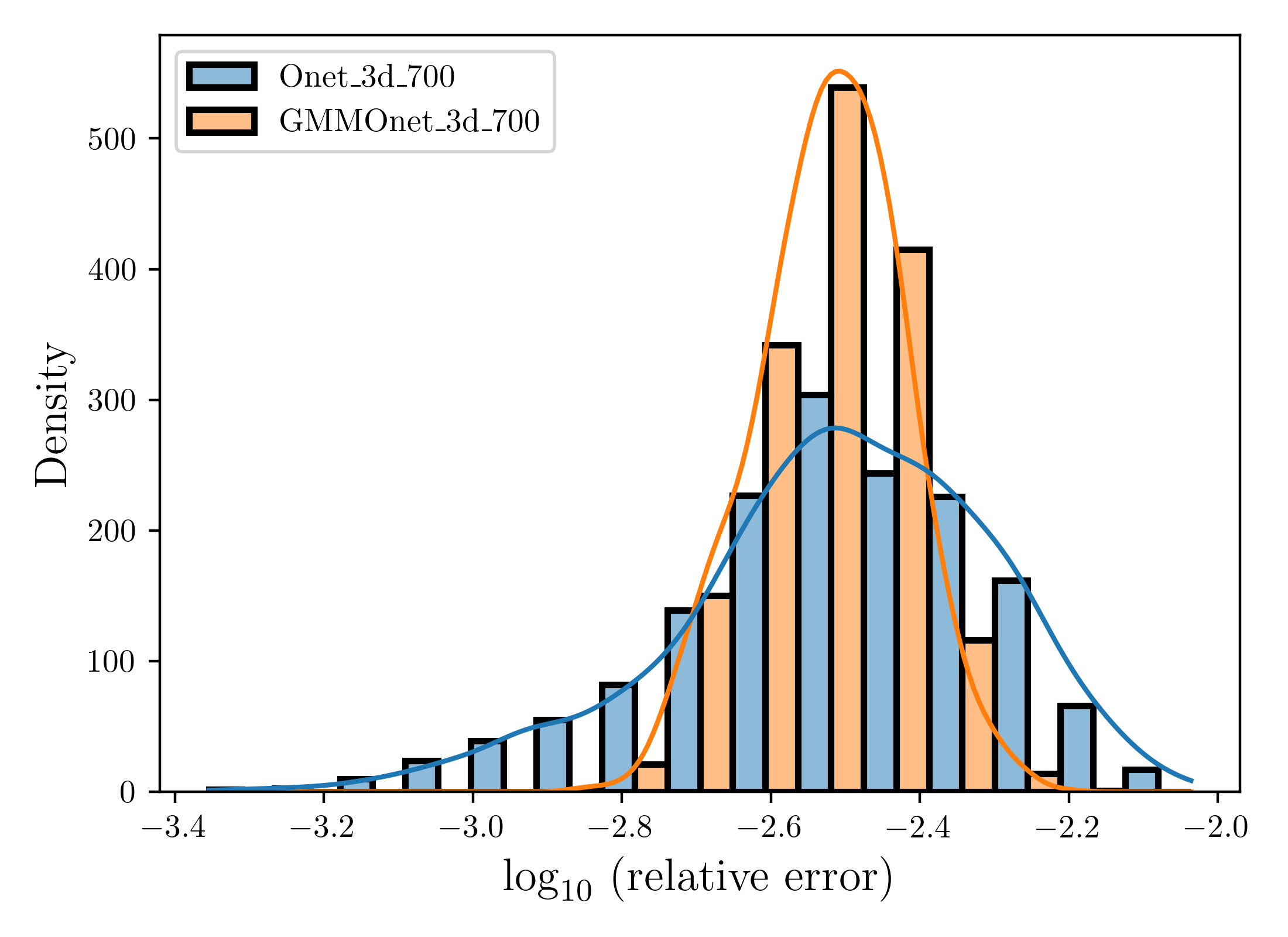}}
    \caption{$\log_{10}(\text{relative error})$ distribution of pressure field predictions at full temporal resolution between model outputs and labeled data. (a),(b),(c): $\log_{10}(\text{relative error})$ distribution on all of the mesh points. (d),(e),(f): $\log_{10}(\text{relative error})$ distribution on the producing well points.}
    \label{fig:3d-re-distribution}
\end{figure}

\begin{figure}[H]
    \centering
    \subfloat[]{\includegraphics[width=0.45\linewidth]{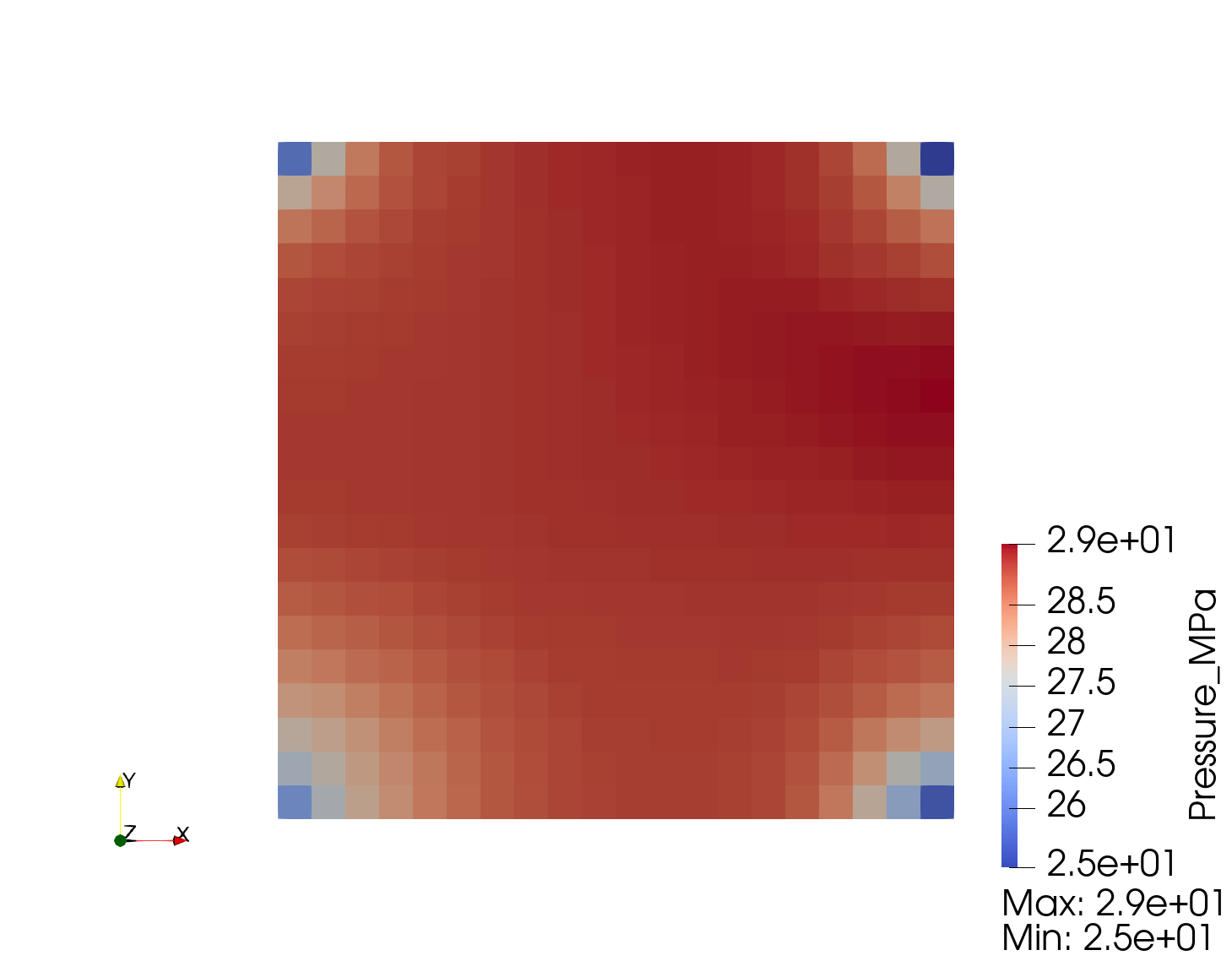}}
    \subfloat[]{\includegraphics[width=0.45\linewidth]{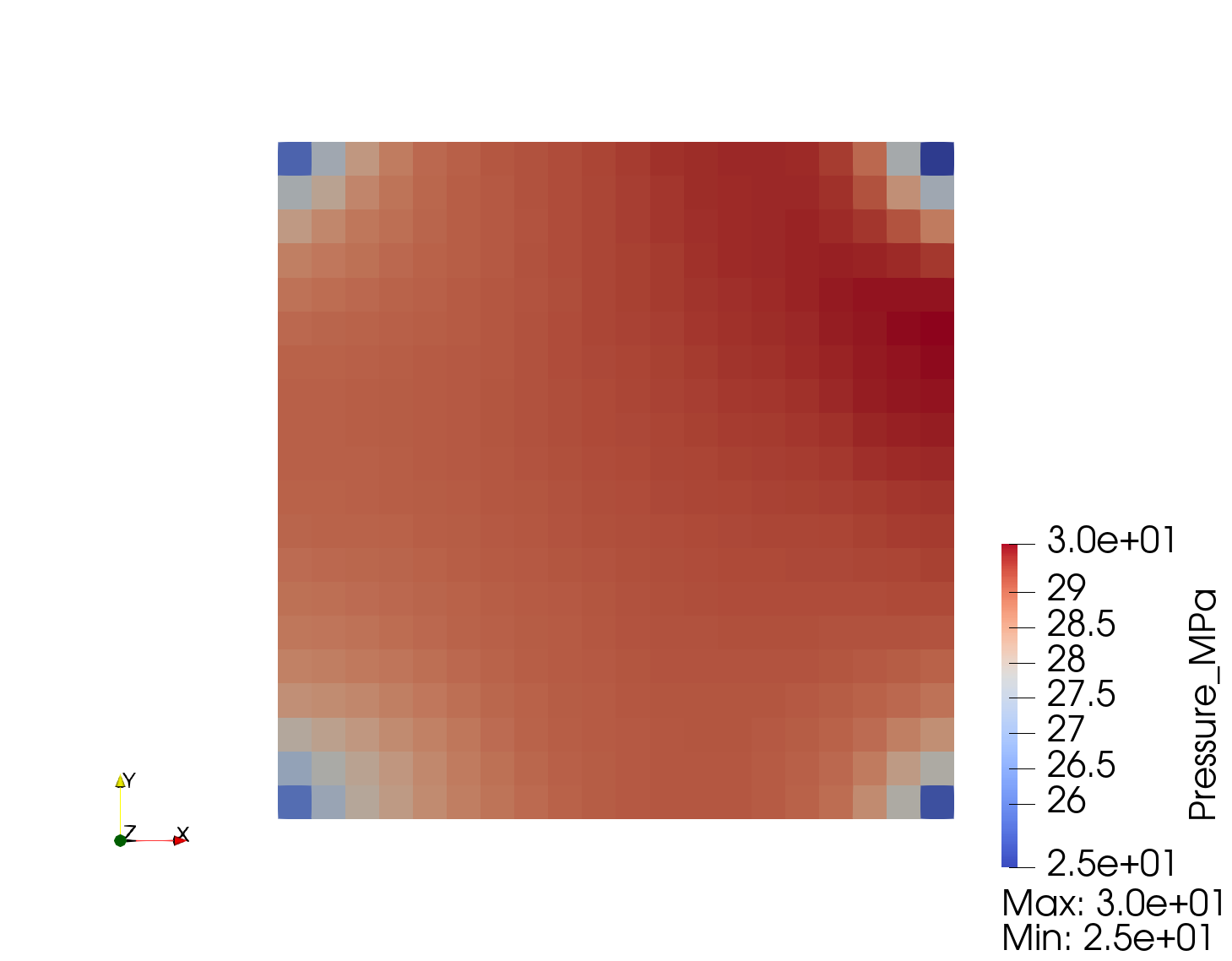}}\\

    \subfloat[]{\includegraphics[width=0.45\linewidth]{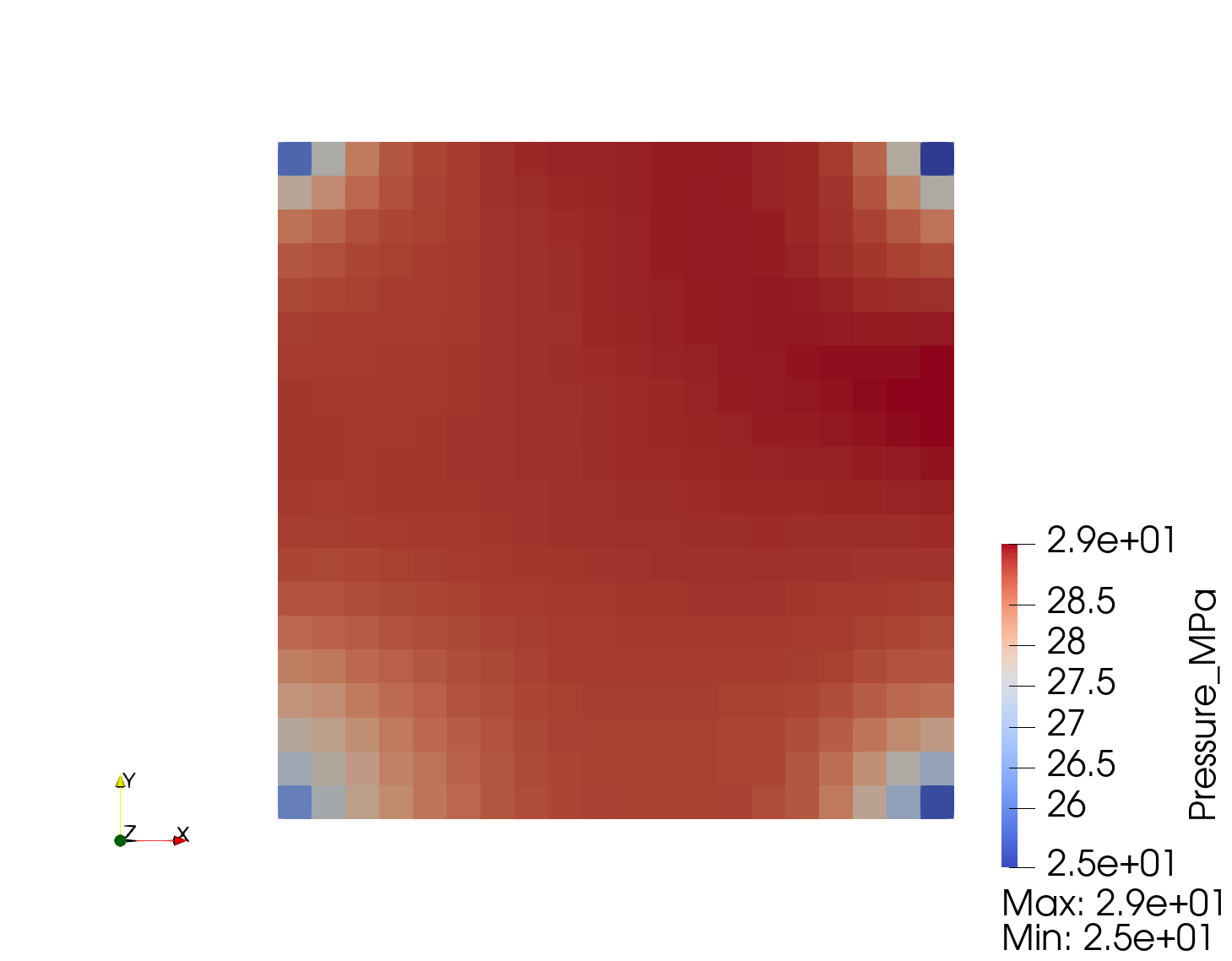}}
    \subfloat[]{\includegraphics[width=0.45\linewidth]{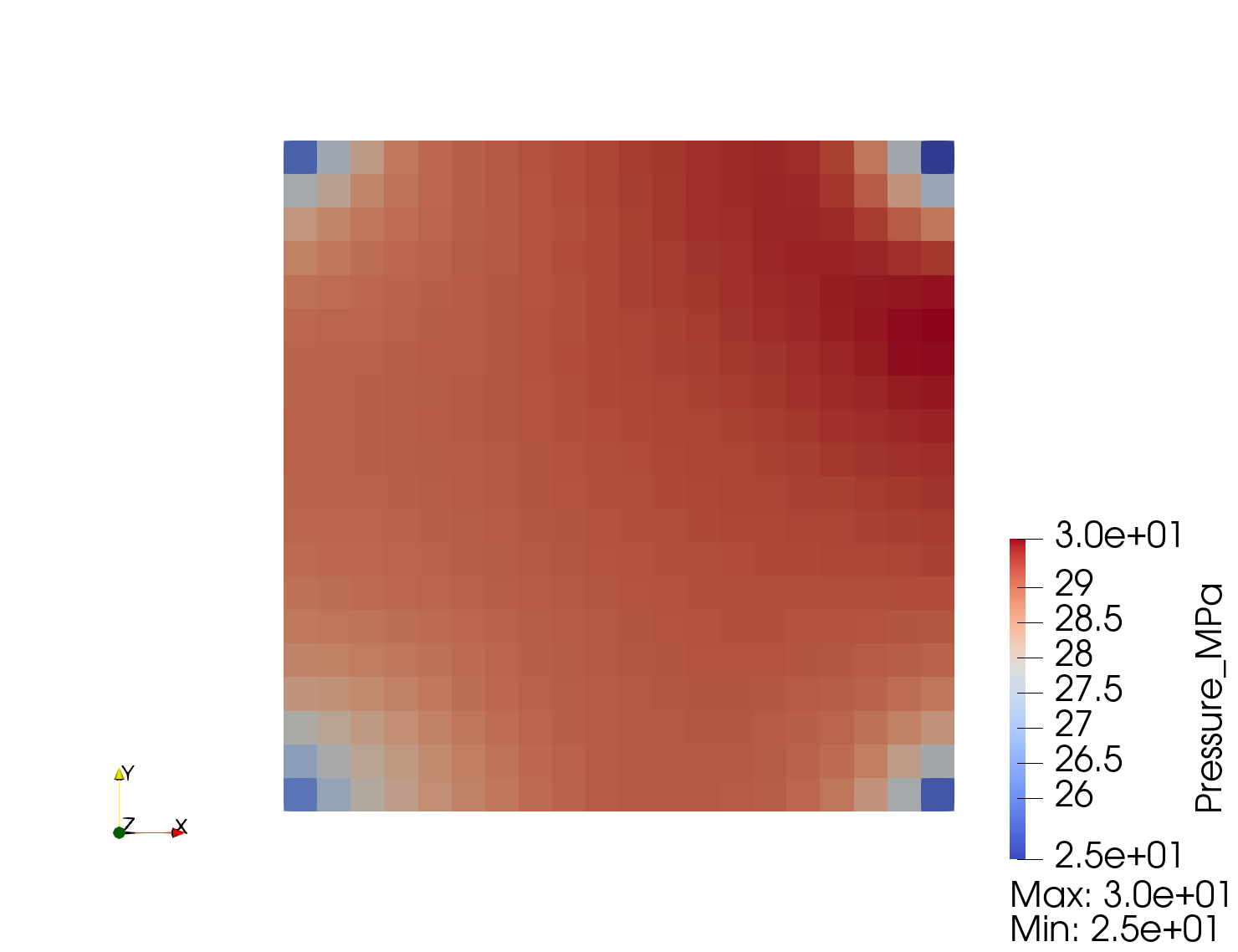}}
    
    \caption{(a),(b): Pressure fields of two distinct 3D permeability field samples at the final time step computed via FVM. (c), (d): Corresponding pressure fields at the final time step computed via AROnet based on adaptive sampling.} 
    \label{fig:3d-predict}
\end{figure}

As demonstrated in \autoref{tab:3d-result} and \autoref{fig:3d-re-distribution}, the surrogate modeling framework with adaptive sampling achieves higher prediction accuracy compared to random sampling consistently. Visual inspections of the predicted pressure fields shows almost no difference between the prediction and ground truth as presented in \autoref{fig:3d-predict}.
    
\subsection{Two-Phase Darcy Flow}\label{sec:2phase-val}
In this test case, there are 4 producing wells where oil is produced over time by the difference of well-block pressure and bottom-hole pressure, and a water-injecting well in the centre of the domain.
The objective is to validate the surrogate for predicting coupled pressure and saturation fields for two-phase Darcy flows. 
The entire dataset comprises 2,100 permeability realizations generated via Gaussian sequential simulation.  Corresponding pressure and water saturation time-series are computed using FVM, with each realization containing 360 temporal snapshots.

\begin{figure}[H]
    \centering
    \subfloat[]{\includegraphics[width=0.33\linewidth]{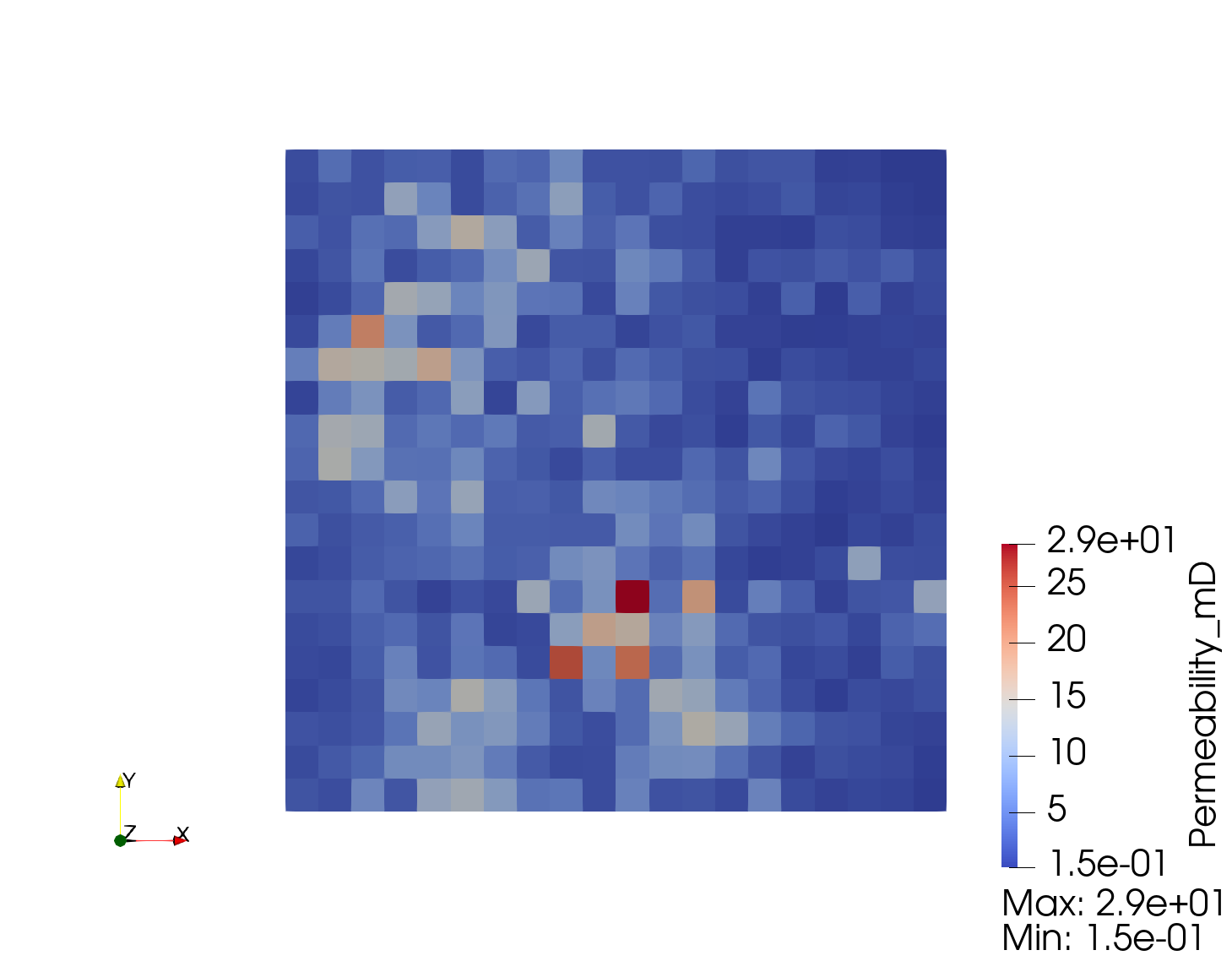}}
    \subfloat[]{\includegraphics[width=0.33\linewidth]{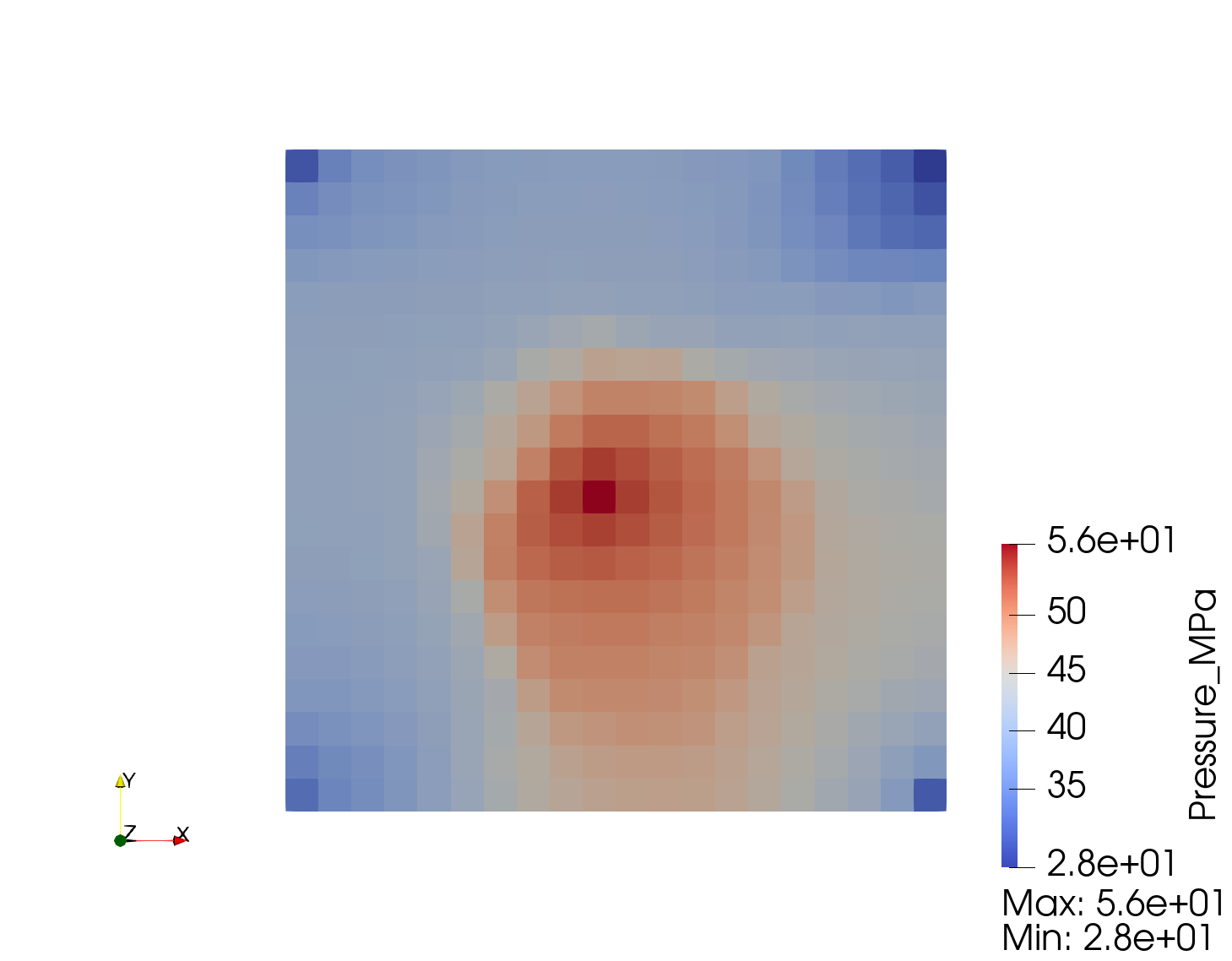}}
    \subfloat[]{\includegraphics[width=0.33\textwidth]{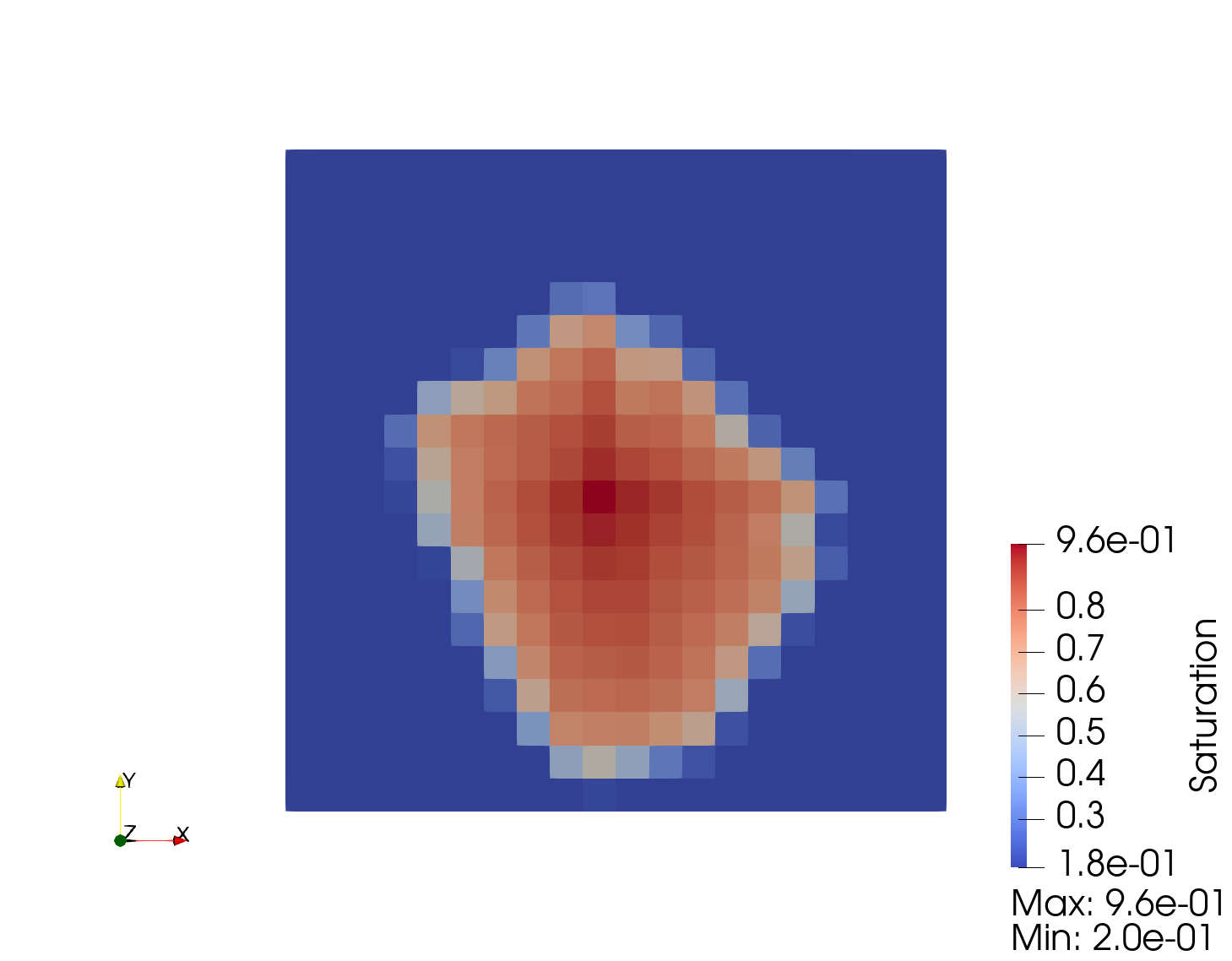}}
    \caption{(a):A sample of two-phase flow's permeability field. (b): corresponding press field computed via FVM at the last time step. (c): corresponding water saturation $S_w$ computed via FVM at the last time step.}
    \label{fig:2phase-perm}
\end{figure}
\begin{figure}[H]
    \centering
    \subfloat[]{\includegraphics[width=0.33\linewidth]{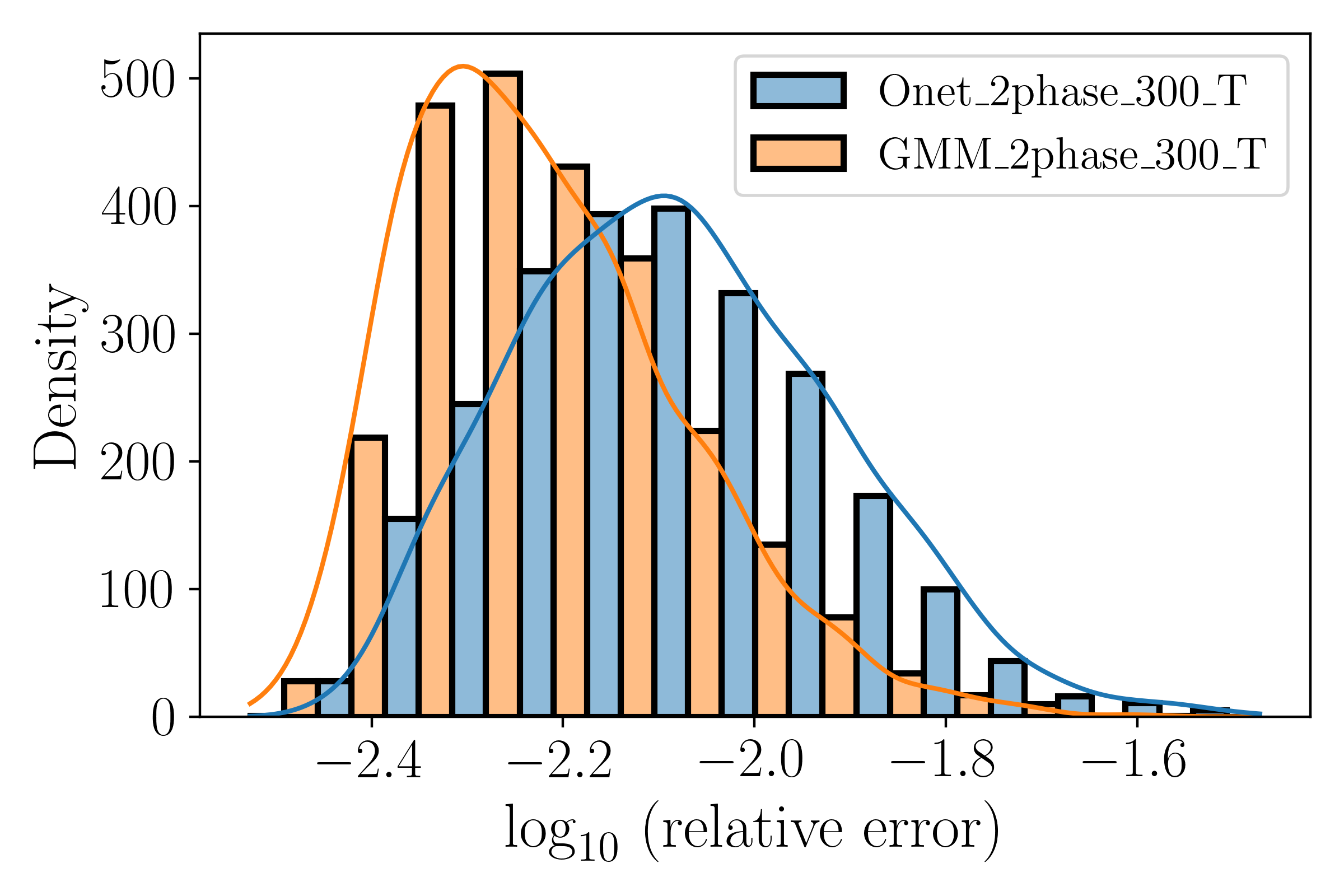}}
    \subfloat[]{\includegraphics[width=0.33\linewidth]{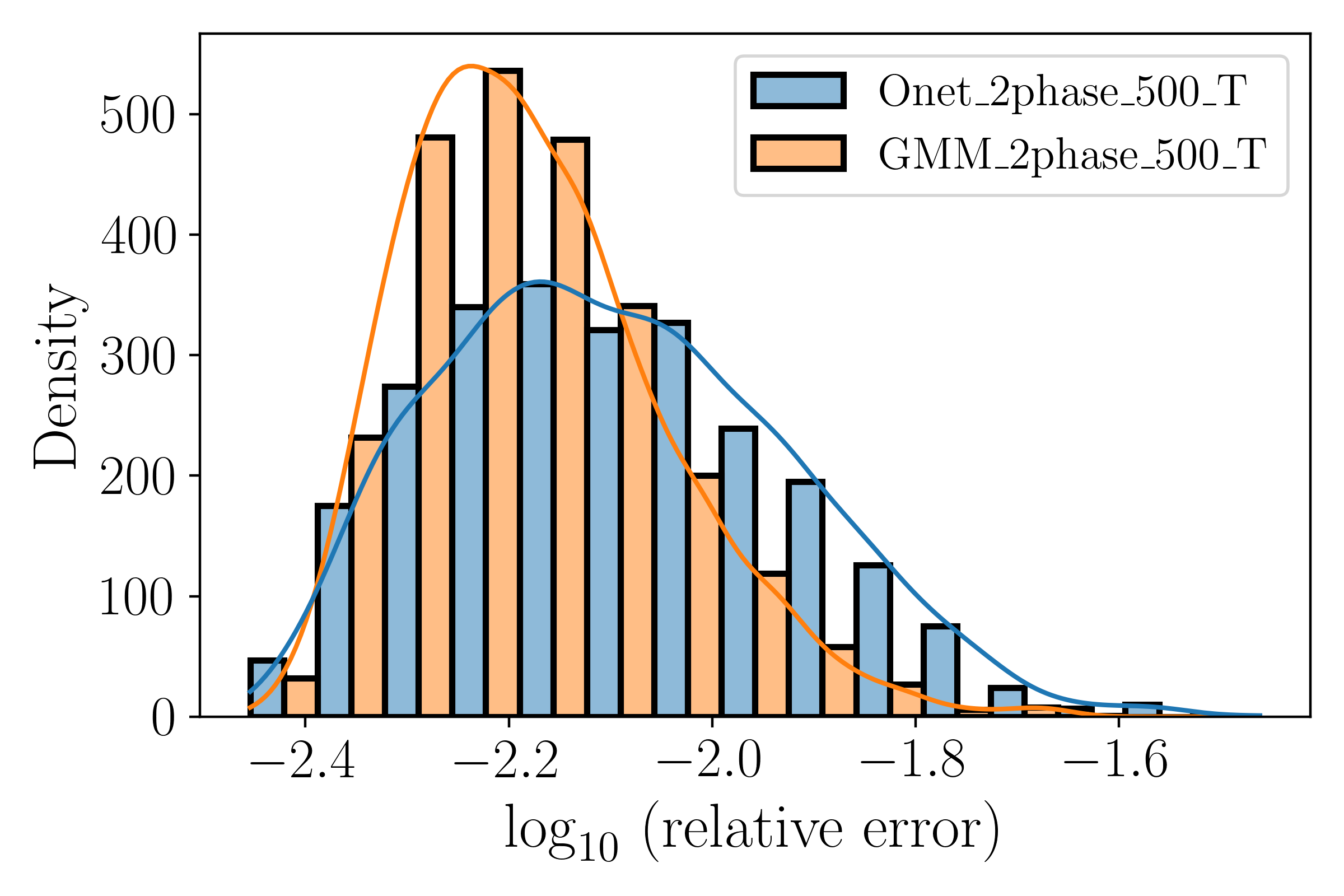}}
    \subfloat[]{\includegraphics[width=0.33\linewidth]{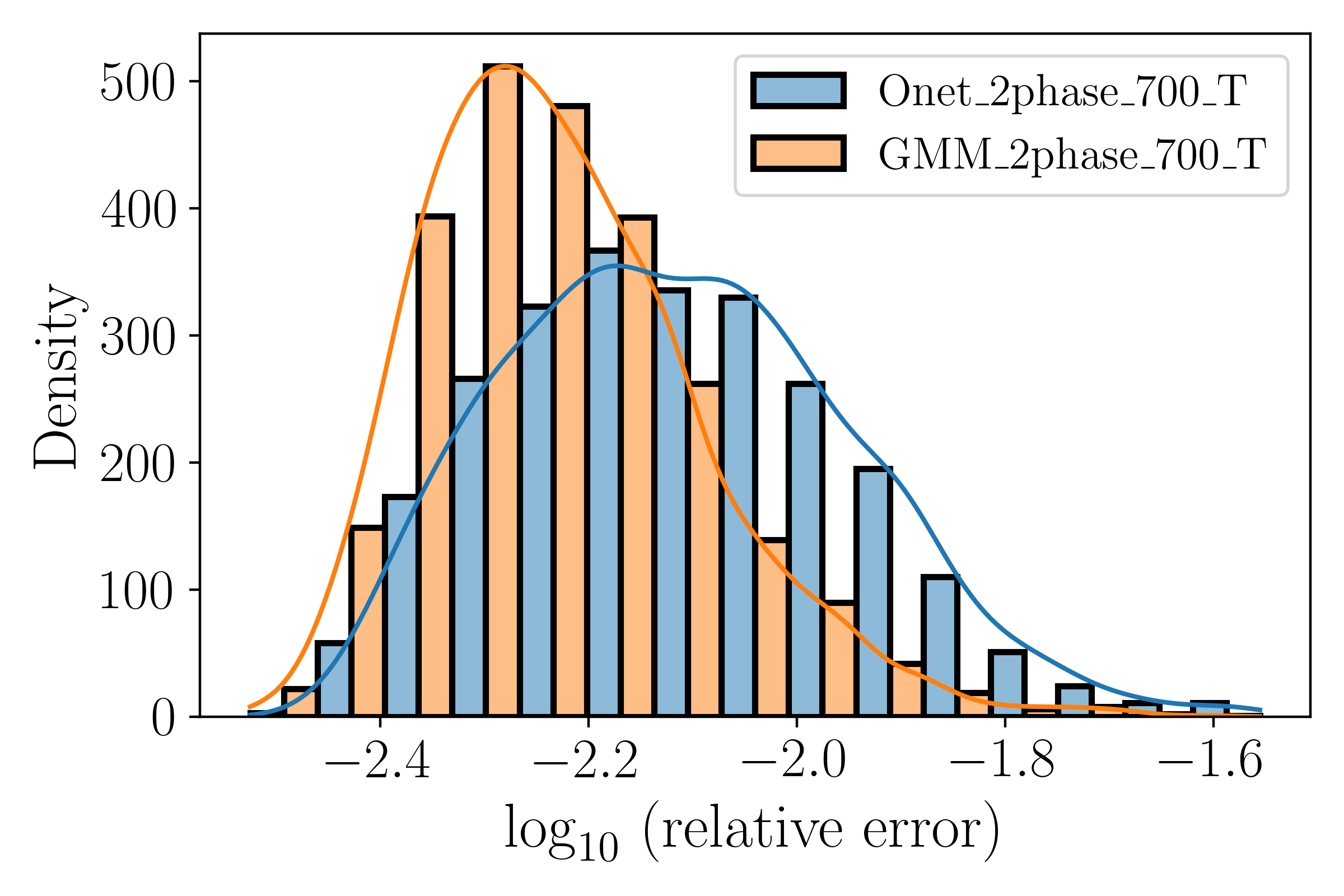}}
    \\
    \subfloat[]{\includegraphics[width=0.33\linewidth]{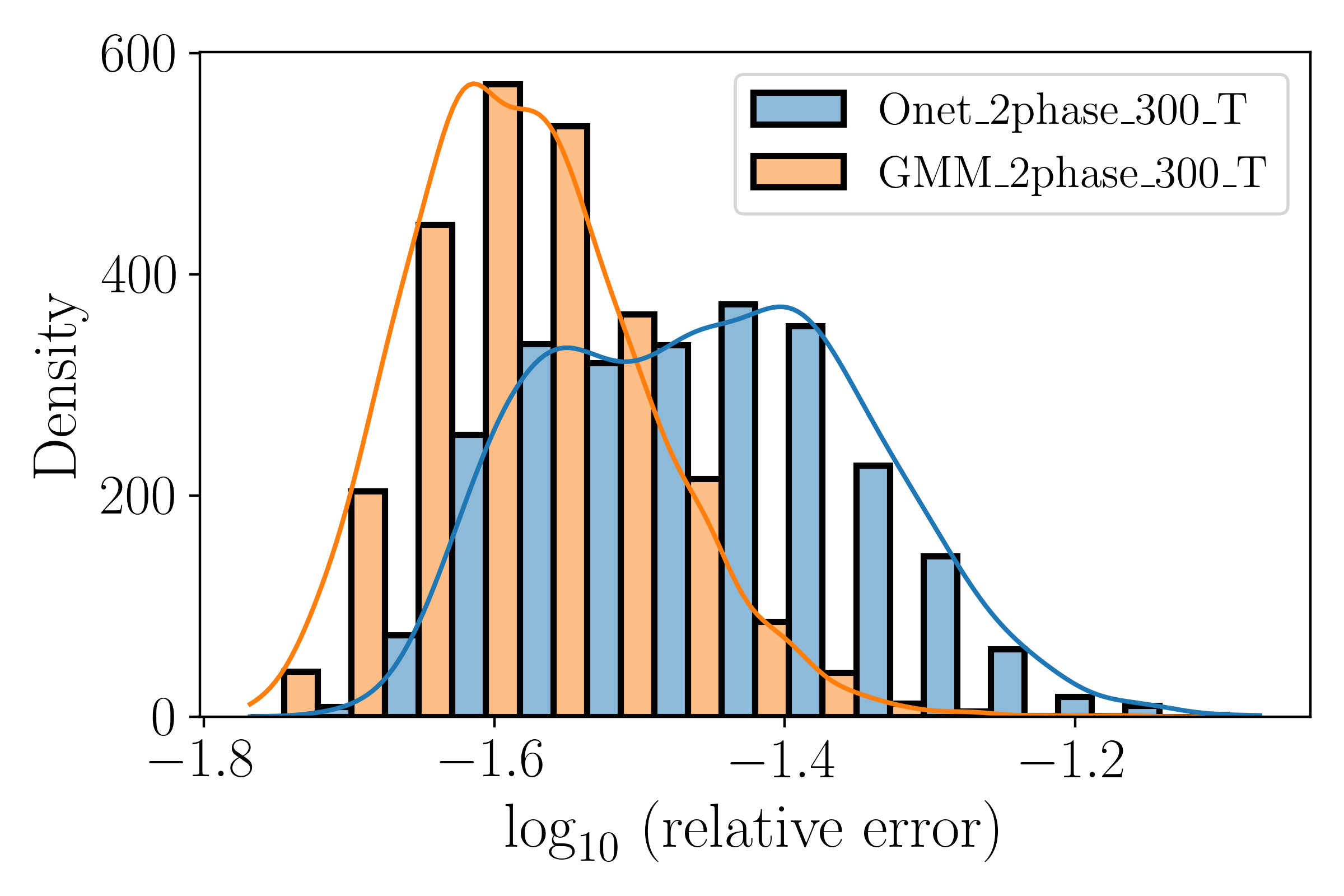}}
    \subfloat[]{\includegraphics[width=0.33\linewidth]{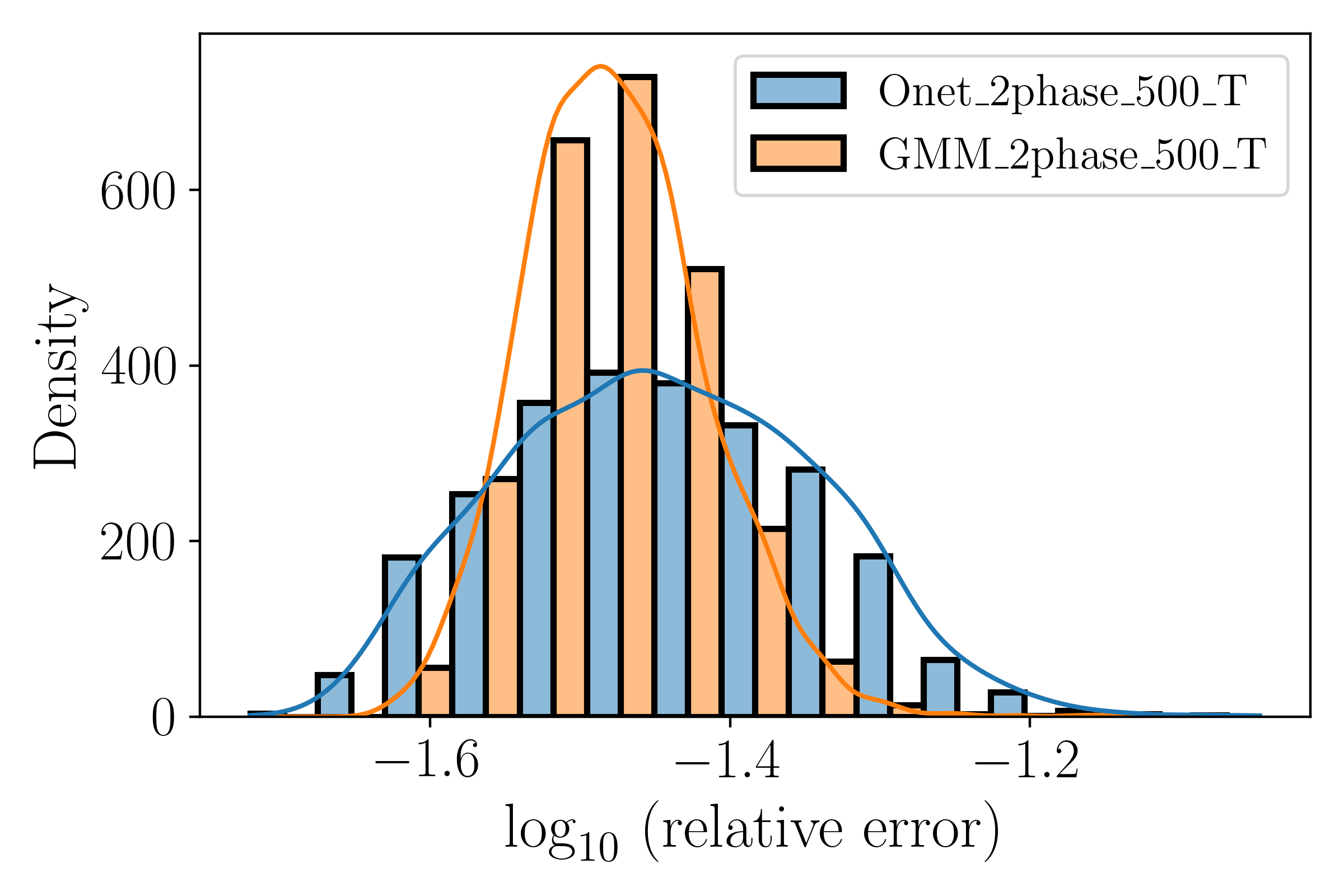}}
    \subfloat[]{\includegraphics[width=0.33\linewidth]{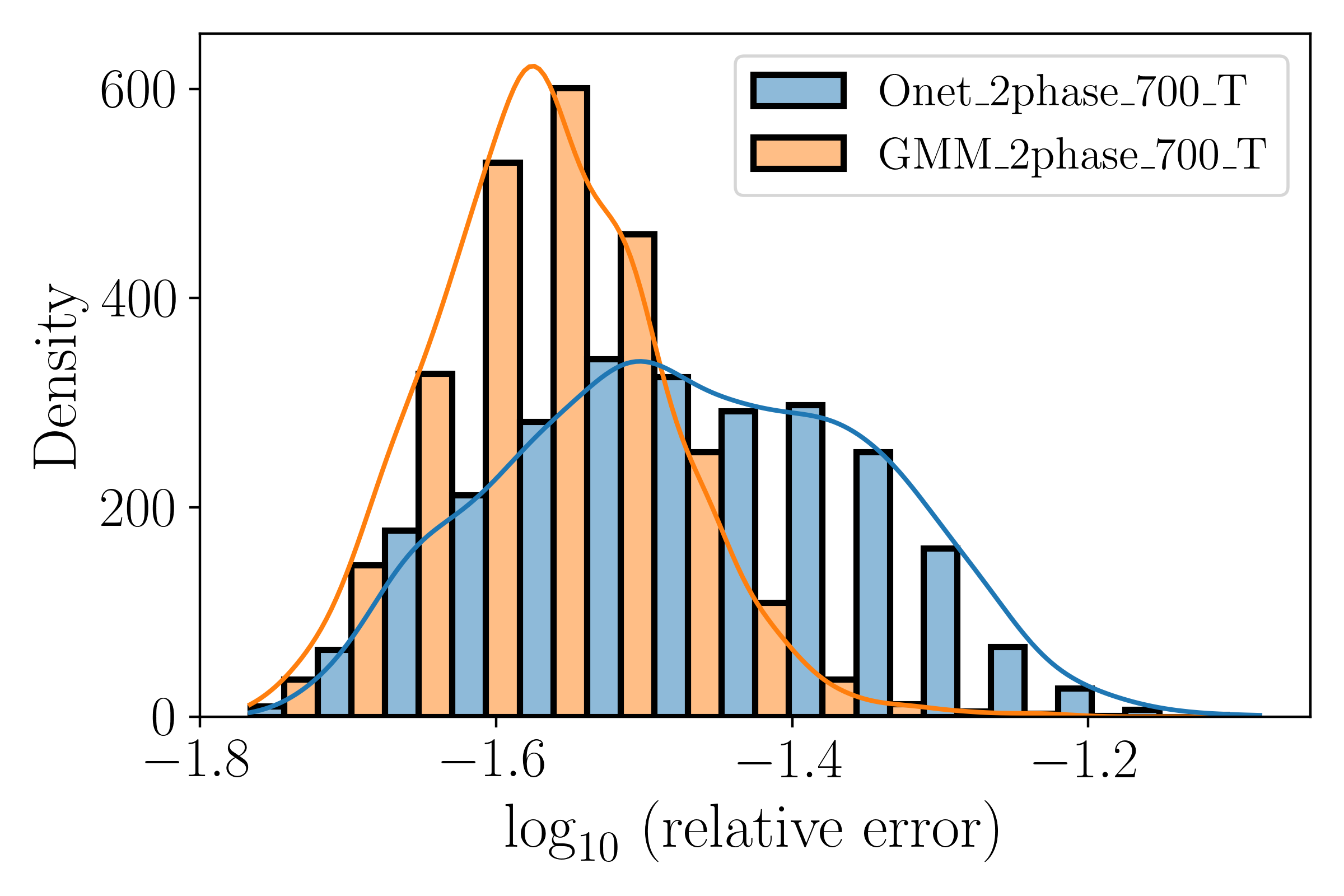}}
    \caption{(a), (b), (c): $\log_{10}(\text{relative error})$ distribution of pressure on all mesh points and the samples of model train set is 300, 500 and 700. (d), (e), (f): $\log_{10}(\text{relative error})$ distribution of water saturation on all mesh points and the samples of model train set is 300, 500 and 700.}
    \label{fig:2phase-result} 
\end{figure}
\begin{figure}[H]
    \centering
    \subfloat[]{\includegraphics[width=0.45\linewidth]{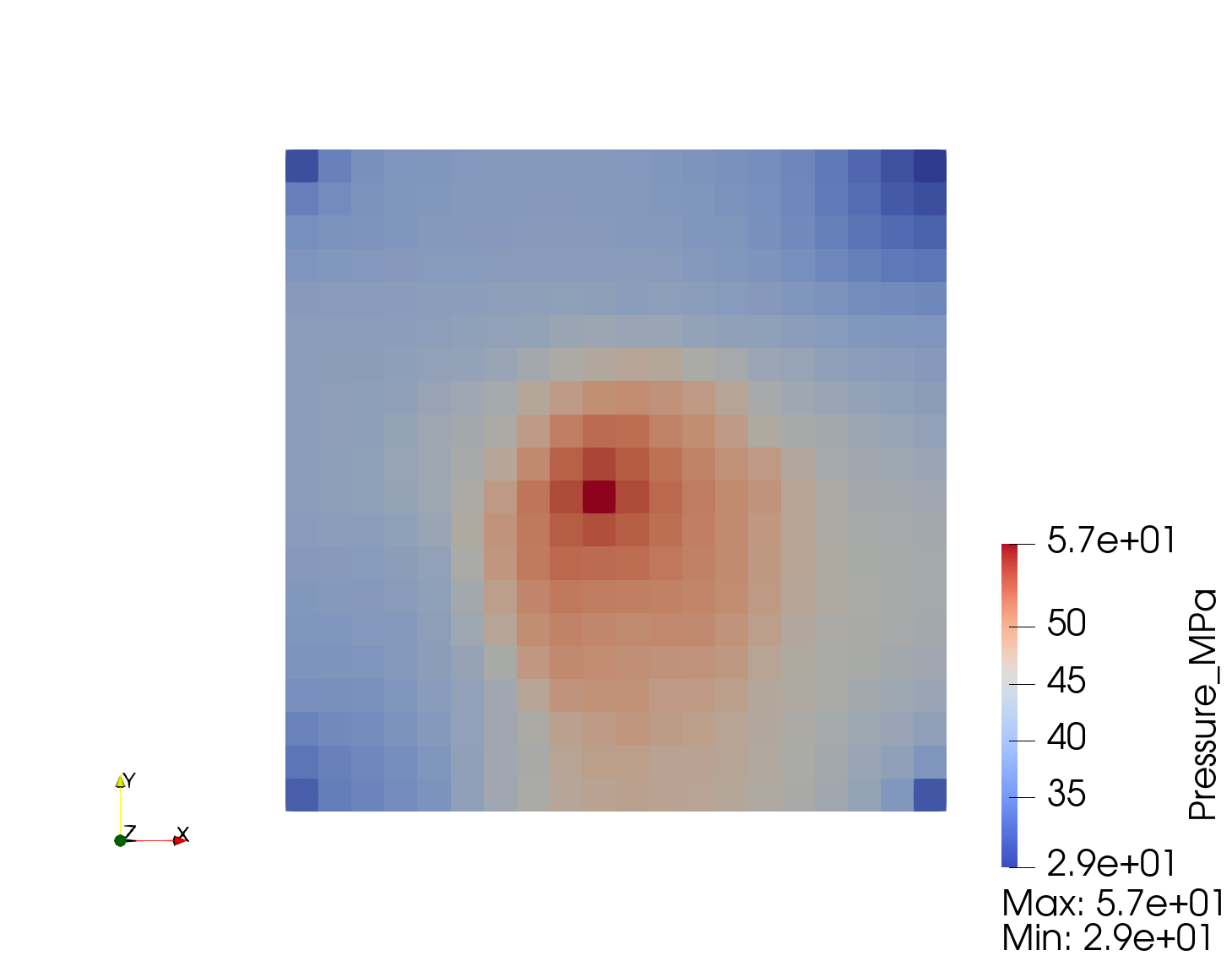}}
    \subfloat[]{\includegraphics[width=0.45\linewidth]{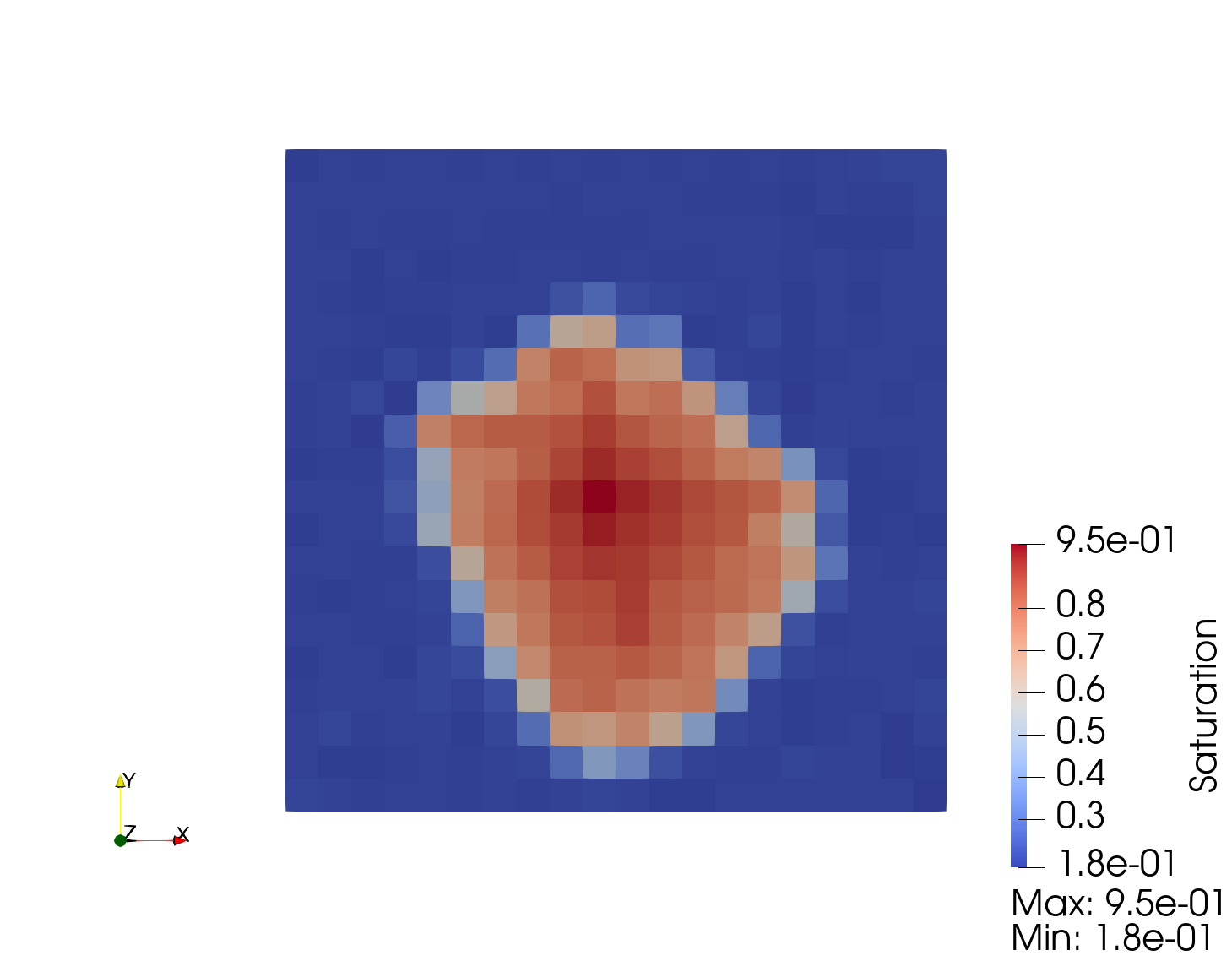}}\\
    \caption{(a): Pressure field predicted by AROnet. (b): Water saturation $S_w$ predicted by AROnet.}
    \label{fig:2phase-Sw_pressure}
\end{figure}

The AROnet with adaptive sampling using GMM is trained to be the surrogate. The model simultaneously predicts pressure and water saturation fields for 15 discrete time steps. The model inputs are transmissibility tensor and time embedding matrices. The outputs are the corresponding press and saturation fields. As demonstrated in \autoref{fig:2phase-perm}, by taking (a) and $t$ as inputs, the outputs (b) and (c) at time $t$ are generated as predicted responses. The training parameters are the same as those specified in \autoref{tab:nn_3d-params} except that Epochs=10. The training set consists of 700 samples with the corresponding pressure and saturation fields for 10 time steps. The predicted results are evaluated in \autoref{tab:2phase-result} and \autoref{fig:2phase-result}, demonstrating the framework's capability for predicting two-phase flow fields. Visual inspections in \autoref{fig:2phase-Sw_pressure} also demonstrates the precision of predicted flow fields.

\begin{table}[H]
    \centering
    \caption{experiment results of two-phase flow test case}
    \begin{tabular}{l l l l}
    \toprule
      samples  &  $S_w$ relative error & press relative error & adaptive sampling \\ 
    \midrule
      300   &  $3.95 \times 10^{-2}$ & $9.01 \times 10^{-3} $ & False \\
      300   &  $3.50 \times 10^{-2}$ & $7.47 \times 10^{-3} $ & True \\
      
      500   &  $3.65 \times 10^{-2}$ & $8.50 \times 10^{-3} $ & False \\
      500   &  $3.41 \times 10^{-2}$ & $7.00 \times 10^{-4} $ & True \\
      
      700   &  $4.10 \times 10^{-3}$ & $8.29 \times 10^{-3} $ & False \\
      700   &  $2.35 \times 10^{-3}$ & $6.25 \times 10^{-3} $ & True \\
    \bottomrule
    \end{tabular}
    \label{tab:2phase-result}
\end{table}

\section{Conclusion}
A discrete neural operator AROnet has been developed for surrogate modeling of spatiotemporal fields in the discrete space associated with numerical solutions of Darcy flows in heterogeneous porous media with random parameters. 
Time embedding is implemented to achieve end-to-end mapping for all time steps. The new discrete neural operator has demonstrated higher prediction accuracy than the SOTA ARUnet structure. Inputs using transmissibility tensors derived from the FVM scheme results in higher accuracy compared to those using the original parameter fields. 
Adaptive sampling has been developed in the latent space to optimize training samples based on estimated probability density estimation for the residuals. Adaptive sampling is realized by training a Gaussian mixture model using expectation maximization to achieve targeted data augmentation, yielding further improvement in prediction accuracy. The algorithm has been validated using 2D/3D single- and two-phase Darcy flows with random heterogeneous parameters.

\section*{CRediT authorship contribution statement}
Zhenglong Chen: Formal analysis, Investigation, Methodology, Validation, Writing-original draft. Zhao Zhang: Conceptualization, Methodology, Validation, Writing – review \&
editing. Xia Yan: Methodology, Validation. Jiayu Zhai: Formal analysis. Piyang Liu: Validation. Kai Zhang: Supervision.

\section*{Declaration of competing interest}
The authors declare that they have no known competing financial interests or personal relationships that could have appeared to influence the work reported in this paper.

\section*{Data availability}
The data that support the findings of this study are available from the corresponding author upon reasonable request.

\section*{Acknowledgements}
The research is supported by the Natural Science Foundation of Shandong Province (No.ZR2024MA057), the Fundamental Research Funds for the Central Universities and the Future Plan for Young Scholars of Shandong University.

\bibliographystyle{elsarticle-harv} 
\bibliography{ref}

@article{huang2021gcaunet,
	title={GCAUNet: A group cross-channel attention residual UNet for slice based brain tumor segmentation},
	author={Huang, Zhicheng and Zhao, Yifan and Liu, Yang and Li, Wei and Zhang, Qiang},
	journal={Biomedical Signal Processing and Control},
	volume={70},
	pages={102958},
	year={2021},
	publisher={Elsevier}
}

@article{ding2024novel,
	title={A novel image denoising algorithm combining attention mechanism and residual UNet network},
	author={Ding, Shuai and Wang, Qian and Guo, Lei and others},
	journal={Knowledge and Information Systems},
	volume={66},
	number={1},
	pages={581--611},
	year={2024},
	publisher={Springer}
}

@article{article,
author = {Karniadakis, George and Kevrekidis, Yannis and Lu, Lu and Perdikaris, Paris and Wang, Sifan and Yang, Liu},
year = {2021},
month = {05},
pages = {1-19},
title = {Physics-informed machine learning},
journal = {Nature Reviews Physics}
}

@article{RAISSI2019686,
title = {Physics-informed neural networks: A deep learning framework for solving forward and inverse problems involving nonlinear partial differential equations},
journal = {Journal of Computational Physics},
volume = {378},
pages = {686-707},
year = {2019},
issn = {0021-9991},
author = {M. Raissi and P. Perdikaris and G.E. Karniadakis},
}

@article{dinh2016density,
  title={Density estimation using real nvp},
  author={Dinh, Laurent and Sohl-Dickstein, Jascha and Bengio, Samy},
  journal={arXiv preprint arXiv:1605.08803},
  year={2016}
}

@article{luLearningNonlinearOperators2021,
  title = {Learning Nonlinear Operators via {{DeepONet}} Based on the Universal Approximation Theorem of Operators},
  author = {Lu, Lu and Jin, Pengzhan and Pang, Guofei and Zhang, Zhongqiang and Karniadakis, George Em},
  year = {2021},
  journal = {Nature Machine Intelligence},
  volume = {3},
  number = {3},
  pages = {218--229},
  publisher = {Nature Publishing Group},
  copyright = {2021 The Author(s), under exclusive licence to Springer Nature Limited}
}

@article{kovachkiNeuralOperatorLearninga,
author = {Kovachki, Nikola and Li, Zongyi and Liu, Burigede and Azizzadenesheli, Kamyar and Bhattacharya, Kaushik and Stuart, Andrew and Anandkumar, Anima},
title = {Neural operator: learning maps between function spaces with applications to PDEs},
year = {2023},
issue_date = {January 2023},
publisher = {JMLR.org},
volume = {24},
number = {1},
issn = {1532-4435},
journal = {J. Mach. Learn. Res.},
month = jan,
articleno = {89},
numpages = {97}
}

@article{wangLearningSolutionOperator2021,
  title = {Learning the Solution Operator of Parametric Partial Differential Equations with Physics-Informed {{DeepONets}}},
  author = {Wang, Sifan and Wang, Hanwen and Perdikaris, Paris},
  year = {2021},
  journal = {Science Advances},
  volume = {7},
  number = {40},
  pages = {eabi8605},
  publisher = {American Association for the Advancement of Science}
}

@article{li2020fourier,
	title={Fourier neural operator for parametric partial differential equations},
	author={Li, Zongyi and Kovachki, Nikola and Azizzadenesheli, Kamyar and Liu, Burigede and Bhattacharya, Kaushik and Stuart, Andrew and Anandkumar, Anima},
	journal={arXiv preprint arXiv:2010.08895},
	year={2020}
}

@article{luDeepXDEDeepLearning2021,
  title = {{{DeepXDE}}: {{A Deep Learning Library}} for {{Solving Differential Equations}}},
  author = {Lu, Lu and Meng, Xuhui and Mao, Zhiping and Karniadakis, George Em},
  year = {2021},
  journal = {SIAM Review},
  volume = {63},
  number = {1},
  pages = {208--228}
}

@article{tangDASPINNsDeepAdaptive2023,
  title = {{{DAS-PINNs}}: {{A}} Deep Adaptive Sampling Method for Solving High-Dimensional Partial Differential Equations},
  author = {Tang, Kejun and Wan, Xiaoliang and Yang, Chao},
  year = {2023},
  journal = {Journal of Computational Physics},
  volume = {476},
  pages = {111868}
}

@article{tangDeepDensityEstimation2020,
  title = {Deep Density Estimation via Invertible Block-Triangular Mapping},
  author = {Tang, Keju and Wan, Xiaoliang and Liao, Qifeng},
  year = {2020},
  journal = {Theoretical and Applied Mechanics Letters},
  volume = {10},
  number = {3},
  pages = {143--148}
}

@article{zhangPhysicsinformedConvolutionalNeural2023,
  title = {A Physics-Informed Convolutional Neural Network for the Simulation and Prediction of Two-Phase {{Darcy}} Flows in Heterogeneous Porous Media},
  author = {Zhang, Zhao and Yan, Xia and Liu, Piyang and Zhang, Kai and Han, Renmin and Wang, Sheng},
  year = {2023},
  journal = {Journal of Computational Physics},
  volume = {477},
  pages = {111919}
}

@article{zhangPhysicsinformedDeepConvolutional2022,
  title = {A Physics-Informed Deep Convolutional Neural Network for Simulating and Predicting Transient {{Darcy}} Flows in Heterogeneous Reservoirs without Labeled Data},
  author = {Zhang, Zhao},
  year = {2022},
  journal = {Journal of Petroleum Science and Engineering},
  volume = {211},
  pages = {110179}
}

@article{tang2023adversarial,
  title={Adversarial adaptive sampling: Unify PINN and optimal transport for the approximation of PDEs},
  author={Tang, Kejun and Zhai, Jiayu and Wan, Xiaoliang and Yang, Chao},
  journal={arXiv preprint arXiv:2305.18702},
  year={2023}
}

@article{oktay2018attention,
  title={Attention u-net: Learning where to look for the pancreas},
  author={Oktay, Ozan and Schlemper, Jo and Folgoc, Loic Le and Lee, Matthew and Heinrich, Mattias and Misawa, Kazunari and Mori, Kensaku and McDonagh, Steven and Hammerla, Nils Y and Kainz, Bernhard and others},
  journal={arXiv preprint arXiv:1804.03999},
  year={2018}
}

@article{cho2014learning,
  title={Learning phrase representations using RNN encoder-decoder for statistical machine translation},
  author={Cho, Kyunghyun and Van Merri{\"e}nboer, Bart and Gulcehre, Caglar and Bahdanau, Dzmitry and Bougares, Fethi and Schwenk, Holger and Bengio, Yoshua},
  journal={arXiv preprint arXiv:1406.1078},
  year={2014}
}

@article{ho2020denoising,
  title={Denoising diffusion probabilistic models},
  author={Ho, Jonathan and Jain, Ajay and Abbeel, Pieter},
  journal={Advances in neural information processing systems},
  volume={33},
  pages={6840--6851},
  year={2020}
}

@article{song2020score,
  title={Score-based generative modeling through stochastic differential equations},
  author={Song, Yang and Sohl-Dickstein, Jascha and Kingma, Diederik P and Kumar, Abhishek and Ermon, Stefano and Poole, Ben},
  journal={arXiv preprint arXiv:2011.13456},
  year={2020}
}

@inproceedings{he2016deep,
  title={Deep residual learning for image recognition},
  author={He, Kaiming and Zhang, Xiangyu and Ren, Shaoqing and Sun, Jian},
  booktitle={Proceedings of the IEEE conference on computer vision and pattern recognition},
  pages={770--778},
  year={2016}
}

@book{scott2015multivariate,
  title={Multivariate density estimation: theory, practice, and visualization},
  author={Scott, David W},
  year={2015},
  publisher={John Wiley \& Sons}
}

@article{peaceman1978interpretation,
  title={Interpretation of well-block pressures in numerical reservoir simulation (includes associated paper 6988)},
  author={Peaceman, Donald W},
  journal={Society of Petroleum Engineers Journal},
  volume={18},
  number={03},
  pages={183--194},
  year={1978},
  publisher={SPE}
}

@article{jiao2024gas,
  title={A Gaussian mixture distribution-based adaptive sampling method for physics-informed neural networks},
  author={Jiao, Yuling and Li, Di and Lu, Xiliang and Yang, Jerry Zhijian and Yuan, Cheng},
  journal={Engineering Applications of Artificial Intelligence},
  volume={135},
  pages={108770},
  year={2024},
  publisher={Elsevier}
}

@book{chen2006computational,
  title={Computational methods for multiphase flows in porous media},
  author={Chen, Zhangxin and Huan, Guanren and Ma, Yuanle},
  year={2006},
  publisher={SIAM}
}

@book{brooks1965hydraulic,
  title={Hydraulic properties of porous media},
  author={Brooks, Royal Harvard},
  year={1965},
  publisher={Colorado State University}
}

@article{TAKBIRIBORUJENI2020104475,
title = {A data-driven surrogate to image-based flow simulations in porous media},
journal = {Computers \& Fluids},
volume = {201},
pages = {104475},
year = {2020},
issn = {0045-7930},
author = {Ali Takbiri-Borujeni and Hadi Kazemi and Nasser Nasrabadi}
}

@article{KAZEMI2023105960,
title = {Physics-informed data-driven model for fluid flow in porous media},
journal = {Computers \& Fluids},
volume = {264},
pages = {105960},
year = {2023},
issn = {0045-7930},
author = {Mohammad Kazemi and Ali Takbiri-Borujeni and Sam Takbiri and Arefeh Kazemi}
}
\newpage
\section{Appendix}

\begin{algorithm}
\SetKwInOut{Input}{input}
\SetKwInOut{Output}{output}
\SetAlgoLined
\Input{samples $Z_j, j=1,\cdots, n$, associated weight $w_j$, numbers of new samples $m$.}
\Output{samples $\hat{Z}_i, i=1, \cdots, m$.}

    compute cumulative distribution function (CDF) using the weights on samples $z_j$\;
    set initial point for resampling $u_1 \sim U(0, \frac{1}{n})$, $c_0=0$\;
    
    \For{$i=1, \cdots, m$}{
        $u_i = u_1 + \frac{i-1}{n}$\;
        \While{$u_i > c_j$}{
            $j = j + 1$\;
        }
        Set $\hat{Z}_i = Z_j$\;
    }
\caption{Systematic resampling}
\label{alg:sys-res}
\end{algorithm}

\begin{algorithm}
\SetKwInOut{Input}{input}
\SetKwInOut{Output}{output}
	
	\Input{initial residual vector $R_0$, initial data vector $Z_0$, numbers of new samples $N$.} 
	\Output{new samples for next training $S$.}
    
	Normalising $R_0$ to be $p, \quad s.t. \quad \int p dx = \sum_i p_i = 1$ \;
    Systematic resampling of $Z_0$ according to the weights $p$, get $\hat{Z}$ \;

    Using $\hat{Z}$ to fit Gaussian mixture model(GMM) with different components and get $p_{\theta_i}, i=1,\cdots, m$ \;
    
    Computing BIC of all $\mathcal{M}_i$ \;
    
    Finding the best model $p_{\theta_i}(x) , i = \arg \min \limits_{i} BIC(\mathcal{M}_i)$ corresponding to the minimum BIC \;

    Sampling N samples $S$ by $p_{\theta}(x)$\;

    \caption{Probability density estimation}
    \label{alg:prob-est} 
\end{algorithm}

\begin{algorithm}
\SetKwInOut{Input}{input}
\SetKwInOut{Output}{output}
\SetAlgoLined

    \Input{initial training data $X_0, y_0$, NN $F_{\theta}$.}
    \Output{trained NN.}
    Set number of training iterations $I$, sampling interval $d$ \;
    Training dataset $X_t = X_0, y_t = y_0$ \;
    \For{$i=1,\cdots,I$ }{
    compute $\Tilde{y} = F(X_t)$, loss function $L=MSE(\Tilde{y} - y)$, residual vector $r^2(x)$ \;
    update $\theta$ by stochastic gradient descent \;
    \If{$i \mod d =0$}{
    Use PCA to transform $X_0$ to $Z_0$ \;
    Use $Z_0$ and $r^2(x)$ to fit adptive sampler $p_{\theta}(x)$ \;
    Use $p_{\theta}(x)$ to sample $m$ samples $S$ \;
    Transform $S$ to $\hat{X}$ to obtain new training dataset $X_t = [X_0; X_1], y_t = [y_0; y_1]$ \;
    }
    }
    \caption{Adaptive sampling in training process}
    \label{alg:as-train} 
\end{algorithm}






\end{document}